\documentclass[11pt, letterpaper, reqno]{amsart}
\usepackage{graphicx, amssymb, color}
\usepackage{amsmath}
\usepackage{latexsym}
\usepackage{hyperref}
\usepackage{dsfont}
\usepackage{amsthm}
\usepackage{tabularx,ragged2e,booktabs,caption}
\usepackage{multirow}

\numberwithin{equation}{section}

\def\tauh{{\hat\tau}}
\def\be{\begin{eqnarray}}
\def\ee{\end{eqnarray}}

\def\b*{\begin{eqnarray*}}
\def\e*{\end{eqnarray*}}

\newtheorem{Theorem}{Theorem}[section]
\newtheorem{Definition}[Theorem]{Definition}
\newtheorem{Proposition}[Theorem]{Proposition}

\newtheorem{Lemma}[Theorem]{Lemma}
\newtheorem{Corollary}[Theorem]{Corollary}
\newtheorem{Remark}[Theorem]{Remark}

\makeatletter \@addtoreset{equation}{section}

\def \E{\mathbb{E}}
\def \F{\mathbb{F}}

\def \L{\mathbb{L}}
\def \M{\mathbb{M}}
\def \N{\mathbb{N}}
\def \P{\mathbb{P}}

\def \R{\mathbb{R}}
\def \S{\mathbb{S}}

\def\Ac{\mathcal{A}}
\def\Bc{\mathcal{B}}

\def\Dc{\mathcal{D}}

\def\Fc{\mathcal{F}}

\def\Jc{\mathcal{J}}

\def\Lc{\mathcal{L}}
\def\Mc{\mathcal{M}}

\def\Oc{\mathcal{O}}

\def\Rc{\mathcal{R}}

\def\ah{{\hat a}}
\def\bh{{\hat b}}
\def\ch{{\hat c}}
\def\fh{{\hat f}}

\def\Xh{{\hat X}}

\def\FF{\mathbf{F}}

\def\TT{\mathbf{T}}

\def\Tr#1{{\rm Tr}\left[#1\right]}

\def \Frac{\displaystyle\frac}

\def\no{\noindent}

\def\={\;=\;}
\def\.{\;.}

\def\vp{\varphi}
\def\eps{\varepsilon}

\def\reff#1{{\rm(\ref{#1})}}

\def\1{\mathds{1}}

\def \ep{\hbox{ }\hfill{ ${\cal t}$~\hspace{-5.1mm}~${\cal u}$   } }

\def \proof{{\noindent \bf Proof. }}
\def \ep{\hbox{ }\hfill$\Box$}

 \def\normeL2#1{\left\|{#1}\right\|_{L^2}}
\newcommand {\lb} {\lambda}
\newcommand {\Chi} {{\bf \raise 1.5pt \hbox{$\delta$}}}

\newcommand{\nb}{\nabla}

\begin{document}

\title[]{A Stochastic Approximation for Fully Nonlinear Free Boundary Parabolic Problems}
\author[]{Erhan Bayraktar}\thanks{E. Bayraktar is supported in part by the National Science Foundation under DMS-0955463,  DMS-1118673, and in part by the Susan M. Smith Professorship.} 
\address[Erhan Bayraktar]{Department of Mathematics, University of Michigan, 530 Church Street, Ann Arbor, MI 48109, USA}
\email{erhan@umich.edu}
\author[]{Arash Fahim}\thanks{A. Fahim is supported by the National Science Foundation under DMS-1209519.}
\address[Arash Fahim]{Department of Mathematics, Florida State University, 1017 Academic Way, Tallahassee, FL 32306, USA}
\email{fahim@math.fsu.edu}


 \maketitle

\begin{abstract}
We present a stochastic numerical method for solving fully non-linear free boundary problems of parabolic type and provide a rate of convergence under reasonable conditions on the non-linearity.
\end{abstract}

\section{Introduction}
\label{Sec:sectint}

 {\color{black} When the option pricing problem is of several dimensions, e.g. basket options,} deterministic methods such as finite difference are almost intractable; because the complexity increases exponentially with the dimension and one almost inevitably needs to use Monte-Carlo simulations. Moreover, many problems in finance, e.g. pricing in incomplete markets and portfolio optimization, lead to fully non-linear PDEs. Only very recently has there been some significant development in numerically solving these non-linear PDEs {\color{black} using Monte Carlo methods}, see e.g. \cite{bt}, \cite{zhang}, \cite{bc}, \cite{mazhang}, \cite{ftw} and \cite{gzz}. When the control problem also contains a stopper, e.g. in determining the super hedging price of an American option, see \cite{MR1809521}, or solving controller-and-stopper games, see \cite{2010arXiv1009.0932B},  the non-linear PDEs have free boundaries.

For solving linear PDEs with free boundaries, i.e. in the problem of American options, Longstaff-Shwartz \cite{longsch}, introduced a stochastic method in which American options are approximated by Bermudan options and least squares approximation is used for doing the backward induction. The major feature in \cite{longsch} is the tractability of the implementation for the scheme proposed in terms of the CPU time in high dimensional problems. The most important feature of this model that facilitates the speed is that the number of paths simulated is fixed. Simulating the paths corresponds to introducing a stochastic mesh for the space dimension and the Bermudan approximation to American options corresponds to time discretization. Stochastic mesh makes sure that  \emph{more important} points in the state space are used in the computation of the value function, an important feature which increases the speed of convergence.
So essentially, this algorithm can be thought of as an explicit finite difference scheme with stochastic mesh. One can in fact prove the convergence rate of the entire ``stochastic" explicit finite difference scheme, see \cite{bw} for a survey of these results and some improvements to the original methodology of Longstaff-Shwartz.

For \emph{semi-linear}  free boundary problems a similar stochastic scheme is given through Reflected Backward Stochastic Differential Equations (RBSDE) in \cite{mazhang} and rate of convergence  is derived to be  $h^{1/4}$ assuming uniform ellipticity for the problem where $h$ the mesh size of the time discretizaton. Here the number of paths, $N$, that one needs to simulate, increases with decreasing $h$ and needs to be chosen in a certain way, see e.g. \eqref{eq:relbhandN}. This is similar to what we have in classical explicit finite difference schemes. To achieve stability, when we decrease the mesh size for time, we need to decrease the mesh size for the space variable to keep the ratio of time step over space step squared in a certain range. As we discussed above, the Monte-Carlo simulation creates a stochastic mesh. The first result in this direction is due to \cite{mazhang}. Later \cite{bc} improved the result of \cite{mazhang} by removing the uniform ellipticity condition. Moreover, they improve the rate of convergence to $h^{1/2}$ by assuming more regularity on the obstacle function.

In this paper, we generalize the Longstaff-Schwartz methodology for numerically solving a large class of \emph{fully non-linear} free boundary problems. We extend the idea in \cite{ftw} to present a stochastic scheme for fully non-linear Cauchy problems with obstacle. As described in Remark \ref{montecarlo}, our scheme is stochastic, i.e. the outcome is a random variable which converges to the true solution pathwise. The convergence of our scheme follows from the methodology of \cite{barlessouganidis}, and the results of \cite{ftw}. For the  convenience of the reader, we sketch the convergence argument in Section \ref{subsec:conv}. 
Under a concavity assumption on the non-linearity and regularity of the coefficients, we obtain a rate of convergence using Krylov's method of shaking coefficients together
with the switching system approximation as in  \cite{bmz}, where a rate of approximation is obtained for classical finite difference schemes for elliptic problems with free boundaries. In \cite{caffarelli-souganidis}, Cafaralli and Souganidis provide a rate of convergence without a concavity assumption on the non-linearity but they consider elliptic problems with non-linearity that depends only on the Hessian.

Appendix \ref{Sec:appendix} is provided to establish the comparison, existence and regularity results for a parabolic switching system with free boundary which is needed to provide the estimations in the rate of convergence proof. This appendix generalizes the result of \cite{jakobsen} for parabolic obstacle problems to parabolic switching systems with obstacle.  Also, it can be considered as the parabolic version of \cite{bmz} where they study elliptic switching systems with  obstacle. Appendix \ref{Sec:appendix_b} contains a proof of the technical Lemma \ref{t-to-T}.

The rest of the paper is organized as follows: In Section \ref{Sec:sectnumericalscheme}, we present the stochastic numerical scheme. In Section \ref{Sec:sectasymptotics}, we present the main results, the convergence rate, and its proof. {\color{black} Section \ref{Sec:numres} is devoted to some numerical results illustrating our theoretical findings.} The appendix  is devoted to the analysis of non-linear switching systems with obstacles, which is an essential ingredient in the proof of our main result and some technical proofs.

\vspace*{0.5cm}

\no {\bf Notation.}\quad 
For scalars $a,b\in\R$, we write $a\wedge b:=\min\{a,b\}$, and $a\vee b:=\max\{a,b\}$. By $\M(n,d)$, we denote the collection of all $n\times d$ matrices with real entries. The collection of all symmetric matrices of size $d$ is denoted $\S_d$, and its subset of nonnegative symmetric matrices is denoted by $\S_d^+$.
For a matrix $A\in\M(n,d)$, we denote by $A^{\rm T}$ its transpose. For $A,B\in\M(n,d)$, we denote $A\cdot B:={\rm Tr}[A^{\rm T}B]$. In particular, for $d=1$, $A$ and $B$ are vectors of $\R^n$ and $A\cdot B$ reduces to the Euclidean scalar product.
For a suitably smooth function $\varphi$ on $Q_T:=(0,T]\times\R^d$, we define 
\b*
|\varphi|_\infty:=\sup_{(t,x)\in Q_T}|\varphi(t,x)|
&\mbox{and}&
|\varphi|_1:=|\varphi|_\infty+\mathop{\sup}\limits_{\tiny
\begin{array}{c}
Q_T\times Q_T\\
x\ne x', t\ne t'
\end{array}
}\frac{|\varphi(t,x)-\varphi(t',x')|}{|x-x'|+|t-t'|^\frac{1}{2}}.
\e*
Finally,  by $\E_{t,x}$ we mean the conditional expectation given $X_t=x$ for a pre-specified diffusion process $X$.

\section{Discretization}
\label{Sec:sectnumericalscheme}
We consider the obstacle problem
\be
\label{equation}
&&\min\left\{-\Lc^X v-F\left(\cdot,v,D v,D^2v\right),v-g\right\}= 0,~~\mbox{on}~[0,T)\times\R^d,\\
&&v=g,~~\mbox{on}~\{T\}\times\R^d, \label{terminal}
\ee
where
\b*
 \Lc^X\varphi
 &:=&
 \frac{\partial\varphi}{\partial t}+\mu\cdot D\varphi+\frac{1}{2}a\cdot D^2\varphi,
 \e*
and
\b*
F:(t,x,r,p,\gamma)\in\R_+\times\Oc\times\R\times\R^d\times\S_d &\longmapsto& F(x,r,p,\gamma)\in\R,
\e*
is a non--linear map, 
$\mu$ and $\sigma$ are maps from $\R_+\times\Oc$ to $\R^d$ and $\M(d,d)$, respectively, $a:=\sigma\sigma^{\rm T}$, $g:[0,T)\times \R^d\to\R$. We consider an $\R^d$-valued Brownian motion $W$ on a filtered probability space $\left(\Omega,\Fc,\F,\P\right)$, where the filtration $\F=\{\Fc_t\}_{t\in[0,T]}$ satisfies the usual conditions, and $\Fc_0$ is trivial. 
For a positive integer $n$, let $h:=T/n$, $t_i=ih$, $i=0,\ldots,n$, and consider the one step Euler discretization
 \begin{equation}\label{Euler}
\Xh_h^{t,x}
 :=
 x+\mu(t,x)h+\sigma(t,x)(W_{t+h}-W_t),
 \end{equation}
of the diffusion $X$ corresponding to the linear operator $\Lc^X$. Then the Euler discretization of the process $X$ is defined by:
\b*
 \Xh_{t_{i+1}}&:=&\Xh^{t_i,\Xh_{t_i}}_{h}.
\e*
We suggest the following approximation of the value function $v$
 \be\label{scheme}
 v^h(T,x):=g(T,x)\!\!
 &\mbox{and}&\!\!
 v^h(t_i,x):=\max\{
\TT_h[v^h](t_i,x),{g}(t_i,x)\}~~\text {for any}~~ x\in\R^d,
\ee
where for a given test function $\psi:\R_+\times\R^d\longrightarrow\R$ we denote
 \be\label{TT}
 \TT_h[\psi](t,x)
 :=
 \E_{t,x}\left[\psi(t+h,\hat X_{t+h})\right] + hF\left(\cdot,\Dc_h\psi\right)(t,x),
 \ee
\be\label{hermit}
 \Dc_h\psi(t_i,x)
 &=&
 \E_{t,x}\left[\psi(t+h,\hat X_{t+h})H_h\right],
 \ee
where $H_h=(H^h_0,H^h_1,H^h_2)^\text{T}$ and
 \b*
 H^h_0=1,
 &H^h_1=\left({\sigma^{\rm T}}\right)^{-1}\;\Frac{W_h}{h},&
 H^h_2=\left({\sigma^{\rm T}}\right)^{-1}\;\frac{W_hW^{\rm T}_h-h\mathbf{I}_d}{h^2}\;\sigma^{-1},
 \e*
provided $\sigma$ is invertible.  {\color{black}Notice that \reff{hermit} comes from
\be\label{int-by-part}
 \E_{t,x}[\psi(N)N]= \E_{t,x}[D\psi(N)]~~~\text{and}~~~
 \E_{t,x}[\psi(N)(N^2-1)]= \E_{t,x}[D^2\psi(N)],
\ee
where $N$ is a standard Gaussian random variable.}
The details can be found in Lemma 2.1 of \cite{ftw}.

\section{Asymptotics of the discrete-time approximation}
\label{Sec:sectasymptotics}
In this section, we present the convergence and the rate of convergence result for the scheme introduced in \reff{scheme}, and the assumptions needed for these results. 
\subsection{The main results}

The proof of the convergence follows the general methodology of Barles and Souganidis \cite{barlessouganidis}, and requires that the nonlinear PDE \reff{equation} satisfies the comparison principle in viscosity sense. 

We recall that an upper-semicontinuous (resp. lower-semicontinuous) function $\underline{v}$ (resp.  $\overline{v}$) on $[0,T]\times\R^d$, is called a viscosity subsolution (resp. supersolution) of \reff{equation} if for any $(t,x)\in[0,T]\times\R^d$ and any smooth function $\varphi$ satisfying
\b*
0=(\underline{v}-\varphi)(t,x)=\max_{[0,T]\times{\R^d}}(\underline{v}-\varphi)\left(\text{resp.}~~0=(\overline{v}-\varphi)(t,x)=\min_{[0,T]\times{\R^d}}(\overline{v}-\psi)\right),
\e*
we have:
\begin{itemize}
\item if $t<T$
\b*
\min\left\{-\Lc^X\varphi-F(\cdot,\Dc \varphi),\varphi-{g}\right\}(t,x)
&\le~\mbox{(resp. $\ge$)}&
0,
\e*
\item if $t=T$, $(\underline{v}-g)(T,x)\le 0$ (resp. $(\overline{v}-g)(T,x)\ge 0$).
\end{itemize}
\begin{Remark}{\rm
Note that the above definition is not symmetric for sub and supersolutions. More precisely, for a subsolution we need to have either 
\b*
-\Lc^X\varphi-F(\cdot,\Dc \varphi)\le0~\text{or}~\varphi-{g}\le0.
\e*
However, for a supersolutions we need to have both
\b*
-\Lc^X\varphi-F(\cdot,\Dc \varphi)\ge0~\text{and}~\varphi-{g}\ge0.
\e*
}
\end{Remark}

\begin{Definition}\label{defcomp}
We say that \reff{equation} has comparison for bounded functions if for any bounded upper semicontinuous subsolution $\underline{v}$ and any bounded lower semicontinuous supersolution $\overline{v}$ on $[0,T)\times{\R^d}$, satisfying
$\underline{v}(T,\cdot)\le \overline{v}(T,\cdot)$,
we have $\underline{v}\le\overline{v}$.
\end{Definition}
We denote by $F_r$, $F_p$ and $F_\gamma$ the partial gradients of $F$ with respect to $r$, $p$ and $\gamma$, respectively. We also denote by $F_\gamma^-$ the pseudo-inverse of the non-negative symmetric matrix $F_\gamma$.

\no {\bf Assumption F}\\
 {\it {\rm (i)}  The nonlinearity $F$ is Lipschitz-continuous with respect to $(x,r,p,\gamma)$ uniformly in $t$, and $|F(\cdot,\cdot,0,0,0)|_\infty<K$ for some positive constant $K$;\\
{\rm (ii)} {\color{black} $\sigma$ is invertible and $|\mu|_1+|\sigma|_1<\infty$.}\\
{\rm (iii)}  
$F$ is elliptic and dominated by the diffusion of the linear operator $\Lc^X$, i.e.
 \be\label{F2}
a^{-1}\cdot\nb_{\!\!\gamma} F\le 1~&\mbox{on}&~\R^d\times\R\times\R^d\times\S_d;
 \ee
{\rm (iv)} $F_p\in{\rm Image}(F_\gamma)$ and $\big|F_p^\text{T}F_\gamma^{-}F_p\big|_\infty<K$;\\
{\rm (v)} 
$ F_r-\Frac{1}{4}F_p^{T}F_\gamma^{-}F_p
 \ge
 0$.
 }\\
\begin{Remark}{\rm 
Assumption {\bf F}(v) is made for the sake of simplicity of the presentation. It implies the monotonicity of the above scheme. If this assumption is not made, one can carry out the analysis in  \cite[Remark 3.13, Theorem 3.12, and Lemma 3.19]{ftw} and approximate the solution of the non-monotone scheme with the solution of an appropriate monotone scheme.
}
\end{Remark}

\begin{Theorem}[Convergence]\label{thmconv}
Suppose that Assumption {\bf F} holds. Also, assume that the fully nonlinear PDE \reff{equation} has comparison for bounded functions. Then for every bounded function ${g}$ Lipschitz on $x$ and $\frac12-$H\"older on $t$, there exists a bounded function $v$ such that
$ v^h \longrightarrow v$ locally uniformly.
Moreover, $v$ is the unique bounded viscosity solution of problem \reff{equation}-\reff{terminal}.
\end{Theorem}
By imposing the following stronger assumption, we are able to derive a rate of convergence for the fully non--linear PDE. 

\no {\bf Assumption HJB}\quad {\it The nonlinearity $F$ satisfies Assumption {\rm\bf F}{\rm (ii)-(v)}, and is of the Hamilton-Jacobi-Bellman type:
 \b*
\frac12 a\cdot\gamma+b\cdot p+ F(t,x,r,p,\gamma)
 &=&
 \inf_{\alpha\in\mathcal{A} }\{\mathcal{L}^{\alpha}(t,x,r,p,\gamma)\},
 \\
 \mathcal{L}^{\alpha}(t,x,r,p,\gamma)
 &:=&
 \frac{1}{2}
 Tr[\sigma^\alpha\sigma^{\alpha{\rm T}}(t,x)\gamma]
 +b^{\alpha}(t,x)p+c^{\alpha}(t,x)r+f^{\alpha}(t,x),
 \e*
where functions $\sigma^\alpha$, $b^\alpha$, $c^\alpha$ and $f^\alpha$ satisfy:
 \b*
\sup_{\alpha\in\Ac}\left(|\sigma^\alpha|_1+|b^\alpha|_1+|c^\alpha|_1+|f^\alpha|_1\right) &<& \infty.
 \e*}
 
\no {\bf Assumption HJB+}\quad {\it The nonlinearity $F$ satisfies {\rm \bf HJB}, and for any $\delta>0$, there exists a finite set $\left\{\alpha_i\right\}_{i=1}^{M_\delta}$ such that for any $\alpha\in\Ac$
 \b*
 \mathop{\inf}\limits_{1\le i\le M_\delta}
 |\sigma^\alpha-\sigma^{\alpha_i}|_\infty+|b^\alpha-b^{\alpha_i}|_\infty+|c^\alpha-c^{\alpha_i}|_\infty+|f^\alpha-f^{\alpha_i}|_\infty
 &\le& \delta.
 \e*
}

\begin{Remark}{\rm
Assumption HJB+ is satisfied if $\Ac$ is a compact separable topological space and $\sigma^\alpha(\cdot)$, $b^\alpha(\cdot)$, $c^\alpha(\cdot)$ and $f^\alpha(\cdot)$  are continuous maps from $\Ac$ to $C_b^{\frac{1}{2},1}$, the space of bounded maps which are Lipschitz in $x$ and $\frac{1}{2}$--H\"older in $t$.
}
\end{Remark}

\begin{Theorem}[Rate of Convergence]\label{thmrateconv}
Assume that the boundary condition $g$ is bounded  Lipschitz on $x$ and $\frac12-$H\"older on $t$. Then, there is a constant $C>0$ such that:
\begin{enumerate} 
\item[{\rm(i)}] under Assumption {\bf\rm HJB}, we have $v-v^h \le Ch^{1/4}$,
\item[{\rm(ii)}] under the stronger condition {\bf\rm HJB+}, we also have $-Ch^{1/10}\le v-v^h$.
\end{enumerate}
\end{Theorem}

It is worth mentioning that in the  finite difference literature, the rate of convergence is usually stated in terms of the discretization in the space variable, i.e. $|\Delta x|$, and the time step, i.e. $|\Delta t|$ equals $|\Delta x|^2$. In our context, the stochastic numerical scheme \reff{scheme} is only discretized in time with time step $h$. Therefore, the rates of convergence in Theorem \ref{thmrateconv} corresponds to the rates $|\Delta x|^{1/2}$ and $|\Delta x|^{1/5}$, respectively.

\subsection{Proof of the convergence result}
\label{subsec:conv}

The proof Theorem \ref{thmconv}, similar to the proof of Theorem 3.6 of \cite{ftw}, is based on the result of \cite{barlessouganidis} which requires the scheme to be consistent, monotone and stable. To be consistence with the notation in \cite{barlessouganidis}, we define
\b*
S_h(t,x,r,\phi):=\min\{h^{-1}(r-\TT_h[\phi](t,x)),r-g(t,x)\},
\e*
and then write the scheme \reff{scheme} as $S_v(t,x,v^h(t,x),v^h)$=0. Notice that by the discussions in \cite{obermanthaleia} and in Section 3 of \cite{oberman}, the consistency and monotonicity for the scheme \reff{scheme} for obstacle problem follows from the consistency and monotonicity of the scheme without obstacle provided by Lemmas 3.11 and 3.12 of \cite{ftw}.
More precisely, we have
\begin{itemize}
\item[(i)] {\bf Consistency.} Let $\varphi$ be a smooth function with bounded derivatives. Then for all $(t,x)\in[0,T]\times\R^d$:
 \begin{align}\label{consistency}
 \lim_{\tiny{\begin{array}{c}
             (t',x')\to(t,x)
             \\
             (h,c)\to (0,0)
             \\
             t'+h\le T
             \end{array}}}
 S_h(t',x',c+\phi(t',x'),c+\phi)
 &=
\min\{ -\left(\mathcal{L}^X\varphi+F(\cdot,\varphi,D\varphi,D^2\varphi)\right),\varphi-g\}(t,x).
 \end{align}
\item[(ii)] {\bf Monotonicity.} Let $\varphi, \psi:~[0,T]\times\R^d\longrightarrow\R$ be two bounded functions. Then:
 \be\label{monotonicity}
 \varphi\le\psi
 &\Longrightarrow&
 S_h(t,x,r,\varphi)\ge S_h(t,x,r,\psi).
 \ee
\end{itemize}
On the other hand, one can show the stability of the scheme \reff{scheme} in the following Lemma. Throughout this section, Assumption of the Theorem \ref{thmconv} are enforced.
\begin{Lemma}\label{vhbdd}
The family $(v^h)_h$ defined by \reff{scheme} is bounded, uniformly in $h$.
\end{Lemma} 
\proof 
Let $C_i=|v^h(t_i,\cdot)|_\infty$. By the argument in the proof of Lemma 3.14 in \cite{ftw}, $|T_h[v^h](t_i,\cdot)|_\infty\le C_{i+1}(1+Ch)+Ch$ {\color{black} where $C>0$ depends only on constant $K$ in assumption {\bf F}}. Therefore,
\b*
C_i\le \max\{|g|_\infty,C_{i+1}(1+Ch)+Ch\}\le \max\{C_{i+1},|g|_\infty\}(1+Ch)+Ch.
\e*
Using a backward induction one could obtain that {\color{black} $C_i\le |g|_\infty e^{CT}+\frac{e^{CT}}{C}$} for some constant $C$ independent of $h$.
\ep\\
The monotonicity, consistency and stability of the scheme result to the following Lemma.
\begin{Lemma}\label{vstar}
Let us define
\b*
\begin{split}
v_*(t,x)&:={\color{black}\lim_{(\delta,h)\to (0,0)}} \inf\{v^{h}(t,y): |x-y|+|s-t| \leq \delta, s \in \{0,h,\cdots\}\cap [0,T] \}, \\
v^{*}(t,x)&:= {\color{black} \lim_{(\delta,h)\to (0,0)}} \sup\{v^{h}(t,y): |x-y|+|s-t| \leq \delta, s \in \{0,h,\cdots\}\cap [0,T] \}.
\end{split}
\e*
Then, $v_*$ and $v^*$ are respectively a viscosity supersolution and a viscosity subsolution of \reff{equation}-\reff{terminal}.
\end{Lemma}
Observe that thanks to Lemma \ref{vhbdd}, $v_*$ and $v^*$ are well-defined and bounded functions and we readily have $v_*\le v^*$. Moreover, functions $v_*$ and $v^*$ are respectively lower semicontinuous and upper semicontinuous. Therefore by the comparison principle for \reff{equation}-\reff{terminal} and Lemma \ref{vstar}, it follows that $v^*=v_*$ and the function $v:=v^*=v_*$ is a viscosity solution of \reff{equation}-\reff{terminal} which completes the proof of Theorem \ref{thmconv}. In addition, uniqueness in the class of bounded functions is a consequence of  comparison principle for the problem.\\
 {\color{black} 
{\bf Proof of Lemma \ref{vstar}.} \\
We show that $v^*$ and $v_*$ are subsolution and super solution at any arbitrary point $(t_0,x_0)\in [0,T]\times\R^d$. We split the proof into the following steps.\\
{\bf Step 1 (${\bf t_0<T}$).} In this case, we only establish the supersolution property of $v_*$. The subsolution property of $v^*$ follows from the same line of arguments. Let $\phi$ be a smooth function such that
\b*
0=\min_{\text{\tiny $[0,T]\times\R^d$}}(v_*-\phi)=(v_*-\phi)(t_0,x_0).
\e*
Since function $v_*$ is bounded, by modifying $\phi$ outside a neighborhood of $(t_0,x_0)$, one can assume that the $(t_0,x_0)$ is a global strict minimum point.  Notice that only the local property of the function $\phi$ matters in the definition of viscosity solution.
Therefore, there exists a sequence $\{(t_n,x_n)\}$, such that $(t_n,x_n)\to(t_0,x_0)$, $v^{h_n}(t_n,x_n)\to v^*(t_0,x_0)$, $\xi_n:=\min(v^{h_n}-\phi)=(v^{h_n}-\phi)(t_n,x_n)\to 0$, and $(t_n,x_n)$ is a global minimum of $v^{h_n}-\phi$. (Obtainig this sequesnce is a standard technique in viscosity solution literature. For more details see \cite{barlessouganidis} and the references therein.) \\
Therefore, $v^{h_n}\ge\phi+\xi_n$. By  the monotonicity of the scheme,  \reff{monotonicity}, we have 
$$S_h(t_n,x_n,v^{h_n}(t_n,x_n),v^{h_n})\le S_h(t_n,x_n,\phi(t_n,x_n)+\xi_n,\phi+\xi_n).$$
 Therefore, by the definition of $v^h$ in \reff{scheme},
\b*
0\le S_h(t_n,x_n,\phi(t_n,x_n)+\xi_n,\phi+\xi_n).
\e*
We divide both sides by $h_n$.  Letting $n \to \infty$ and using \reff{consistency} we obtain: 
\b*
0\le\min\{-\left(\mathcal{L}^X\varphi+F(\cdot,\varphi,D\varphi,D^2\varphi)\right),\varphi-g\}(t_0,x_0).
\e*
{\bf Step 2 (Supersolution property when ${\bf t_0=T}$).} Observe that since $g\le v^h$, we have
\b*
g(T,x_0)\le\mathop{\liminf}\limits_{
             (h,t',x')\to(0,T,x_0)
             }
 v^h(t',x')=v_*(T,x_0),
\e*
which completes the proof of the supersolution argument at terminal time.\\
{\bf Step 3 (Subsolution property when ${\bf t_0=T}$).}
Observe that by the definition of $v^h$ and $v ^*$, we readily have $v^*\ge g$. To complete the subsolution argument, we have to show that $v^*(T,\cdot)= g(T,\cdot)$. It is sufficient to show that 
\begin{Lemma}\label{t-to-T}
For all $x\in \R^d$ and $i=0,\cdots,n$, we have
\b*
|v^h(t_i,x)-g(T,x)|\le C\sqrt{T-t_i}.
\e*
\end{Lemma}
The proof of Lemma \ref{t-to-T} is similar to Lemma 4.2 \cite{fahim} or Lemma 3.17 in \cite{ftw}. For the convenience of the reader, we adjust the proof for free boundary problems. However, because the proof is technical and not related to the main result of the paper, we prefer to present it in Appendix \ref{Sec:appendix_b}.\\


\subsection{Derivation of the rate of convergence}
\label{proofrate}

The proof of Theorem \ref{thmrateconv} is based on Barles and Jakobsen \cite{barlesjakobsen}, which uses switching systems approximation and the Krylov method of shaking coefficients \cite{krylov}, \cite{krylov0}, \cite{krylov2}, and  \cite{krylov1}. This has been adapted to classical finite difference schemes for elliptic obstacle (free boundary) problems in \cite{bmz}.  In order to use the method, we need to introduce a comparison principle for the scheme which we will undertake next.
\begin{Proposition} \label{propmaxpri}
Let Assumption $\FF$ hold and set $\beta:=|F_r|_\infty$. Consider two arbitrary bounded functions $\varphi$ and $\psi$ satisfying:
 \be\label{hyplemmaxpri}
\min\left\{ h^{-1}\left(\varphi-{\TT}_h[\varphi]\right),\vp-g\right\} \;\le\; g_1
 &\mbox{and}&
 \min\left\{ h^{-1}\left(\psi-{\TT}_h[\psi]\right) ,\psi-g\right\}\;\ge\; g_2
 \ee
for some bounded functions $g_1$ and $g_2$. Then, for every $i=0,\cdots,n$:
 \be\label{conlemmaxpri}
 (\varphi-\psi)(t_i,x)
 &\le&
 e^{\beta(T-t_i)}\left(|(\varphi-\psi)^+(T,\cdot)|_{\infty}
 +(1+T-t_i)|(g_1-g_2)^+|_\infty\right).
 \ee
\end{Proposition}

\proof 
This follows from Lemma 2.4 of \cite{jakobsen} where it is provided for general monotone schemes for obstacle problems. 

\subsubsection{Proof of Theorem \ref{thmrateconv} (i)}
Under  Assumption HJB, we can build {\color{black} a bounded subsolution ${v}^\eps$ of the nonlinear PDE}, by the method of shaking coefficients, see \cite{barlesjakobsen}, \cite{bmz}, \cite{krylov1}, and the references therein.  More precisely, consider the following equations
\be\label{shake}
\min\left\{-\Lc^X v-\inf_{\tiny
0<s<\eps^2,
|y|<\eps
}\!\!\!\!\!F\left(t-s,x+y,v,D v,D^2v\right),v-g\right\}&= 0,~~\mbox{on}~[0,T)\times\R^d,\\
v&=g,~~\mbox{on}~\{T\}\times\R^d\label{shaketerminal}.
\ee
By Theorem \ref{switchexist}, there exists a unique bounded solution ${v}^\eps$ to \reff{shake}-\reff{shaketerminal}.\\ {\color{black}Because $\inf_{\tiny
0<s<\eps^2,
|y|<\eps
} F \le F$, ${v}^\eps$ is a subsolution to \reff{equation}-\reff{terminal} at all $(t,x)\in[0,T]\times\R^d$ and by Theorem \ref{contdep} (continuous dependence for HJB equations), $v^\eps$ approximates $v$ uniformly, i.e., there exists a positive constant $C$ such that   $v-C\eps \le v^\eps \le v.$}\\
Let $\rho(t,x)$ be a $C^\infty$ non-negative  function supported in $\{(t,x):t\in[0,1],|x|\le 1\}$ with unit mass, and define
 \be\label{defoverueps}
{w}^\eps(t,x)
 :=v^\eps*\rho^\eps
 &\mbox{where}&
 \rho^\eps(t,x)
 :=\frac{1}{\eps^{d+2}}\rho\left(\frac{t}{\eps^2},\frac{x}{\eps}\right).
 \ee 
It follows that $|{w}^\eps-v|\le C\eps$. From the concavity of the the nonlinearity $F$, and Lemma A.3  in \cite{barlesjakobsen2}, ${w}^\eps\in C^\infty$, ${w}^\eps$ is a classical subsolution of \reff{equation} on $U:=\{~(t,x)~|~g(t-s,x+y)<v^\eps(t-s,x+y);\text{   for any  } s\in[0,\eps^2)\text{ and }|y|<\eps\}$. 
\footnote{This heuristically follows from  $
F\left(\cdot,\cdot,v^\eps,D v^\eps,D^2v^\eps\right)*\rho_\eps (t,x)\le F\left(t,x,w^\eps,D w^\eps,D^2w^\eps\right)
$.}
Moreover, by Theorem \ref{switchreg}, $v^\eps$ is Lipschitz in $x$ and $1/2-$H\"older continuous in $t$.  Thus, by Theorem 2.1 in \cite{krylov2}, 
 \be\label{regunderweps}
  \left|\partial_t^{\beta_0}D^\beta{w}^\eps\right| 
 \le C\eps^{1-2\beta_0-|\beta|_1}
 &\mbox{for any}&(\beta_0,\beta)\in\N\times\N^d\setminus\{0\},
 \ee
where $|\beta|_1:=\sum_{i=1}^d\beta_i$, and $C>0$ is some constant. 
As a consequence of the consistency of $\TT_h$, see Lemma 3.22 of \cite{ftw}, we know that
 \b*
 \Rc_h[{w}^\eps](t,x)
 &:=&
 \frac{{w}^\eps(t,x)-\TT_h[{w}^\eps](t,x)}{h}
 \;+\;
 \mathcal{L}^X{w}^\eps(t,x)+F(\cdot,{w}^\eps,D{w}^\eps,D^2{w}^\eps)(t,x)\le  Ch\eps^{-3}.
 \e*
From this estimate together with the subsolution property of $w^\eps$, we see that ${w}^\eps \le {\TT}_h[{w}^\eps]+ Ch^2\eps^{-3}$ holds true on $U$. In addition, by the regularity properties of $g$, one can see that ${w}^\eps\le g+C\eps$ on $[0,T]\times\R^d\setminus U$. Therefore, 
\b*
\min\left\{\frac{{w}^\eps(t,x)-\TT_h[{w}^\eps](t,x)}{h},{w}^\eps-g\right\}\le C_1(\eps+\eps^{-3}h).
\e*
Then, it follows from Proposition \ref{propmaxpri} that
 \be\label{underlineww}
 {w}^\eps-v^h\le C|({w}^\eps-v^h)(T,.)|_\infty+C_1(\eps+h\eps^{-3})\le C(\eps+h\eps^{-3}).
 \ee
Therefore, $ v-v^h
 \le
 v-{w}^\eps+{w}^\eps-v^h
 \le C(\eps + h\eps^{-3})$.
Minimizing the right hand-side estimate over $\eps>0$, we obtain $v-v^h \le Ch^{1/4}$.\ep

\subsubsection{Proof of Theorem \ref{thmrateconv} (ii)}
\label{subsec:thmrateconvii}

To prove the lower bound on the rate of convergence, we will use Assumption {\bf HJB+} and build a switching system approximation to the solution of the nonlinear obstacle problem \eqref{equation}-\eqref{terminal}. This proof method has been used for Cauchy problems of \cite{barlesjakobsen} and \cite{ftw}. For obstacle problems, this method is used in the elliptic case by \cite{bmz} for the classical finite difference schemes. We apply this methodology for parabolic obstacle  problems to prove the lower bound for the convergence rate of our stochastic finite difference scheme.  We split the proof into the following steps:
\begin{enumerate}
\item Approximating the solution to \reff{equation}-\reff{terminal} by a switching system, which relies on Theorem~\ref{contdep}, the continuous dependence result for switching systems with obstacle.
\item Building an almost everywhere smooth supersolution to  \eqref{equation}-\eqref{terminal}  using the mollification of the solution to the switching system. 
\item Using Proposition~\ref{propmaxpri}, the comparison principle for the scheme, to bound the difference of the supersolution obtained in Step 2 and the approximate solution obtained from the scheme.
\end{enumerate}

\no{\bf Step 1.}  Consider the following switching system:
\begin{align}\label{switchshaken}
\min\left\{\max\left\{\!-v^{\eps,i}_t{\color{black}+}\sup_{0<s<\eps^2,
|y|<\eps
}\Lc^{\alpha_i} v^{\eps,i}(\cdot\!-\!s,\cdot\!+\!y),v^{\eps,i}\!\!-\!\!\Mc^{(i)}v^\eps\right\},v^{\eps,i}-g\right\}(t,x)&=0,\\
v^{\eps,i}(T,\cdot)&=g(T,\cdot),\label{switchshakenterminal}
\end{align}
where $v^\eps=(v^{\eps,i})_{i=1}^M$, $\Mc^{(i)}v^\eps=\min_{j:j\ne i}\{v^{\eps,j}+k\}$, $k$ is a non-negative constant, $\alpha_i$'s,  for $i=1,\cdots,M$, are as in assumption {\bf HJB+}, and  $\Lc^{\alpha_i}\vp:=\frac12\Tr{a^{\alpha_i}(t,x)D^2\vp}+b^{\alpha_i}(t,x)D\vp+ c^{\alpha_i}(t,x)\vp+f^{\alpha_i}(t,x)$. 

{\color{black} The above system of equations approximates \reff{equation}-\reff{terminal}. Intuitively speaking,
Assumption {\bf HJB+} introduces a set of approximating controls $\{\alpha_i\}_{i=1}^{M_\delta}$ in $\Ac$. In the corresponding optimization problem, the maximum cost of restricting controls to the set $\{\alpha_i\}_{i=1}^{M_\delta}$ is proportional to $\delta$. In addition, the above switching system imposes a switching cost of $k$ between controls in the finite set $\{\alpha_i\}_{i=1}^{M_\delta}$. If $k$ goes to zero, then all functions $v^{\eps,i}$ in the solution of problem \reff{switchshaken}-\reff{switchshakenterminal} converges to the function $v^\eps$, the solution of the problem without switching cost, i.e. \reff{shake}-\reff{shaketerminal}. On the other hand, we have already seen that  $v^\eps$ approximates function $v$, the solution of \reff{equation}-\reff{terminal}.
}

More rigorously, by Theorem \ref{switchexist} the viscosity solution $\left(v^{\eps,i}\right)_{i=1}^M$  to \reff{switchshaken}-\reff{switchshakenterminal} exists and by Theorem \ref{switchreg} is Lipschitz continuous on $x$ and $\frac12$-H\"older continuous on $t$. Moreover, by using Assumption {\bf HJB+}, Theorem \ref{contdep} and Remark \ref{infsup}, one can approximate the solution to \reff{equation}-\reff{terminal} by the solution to \reff{switchshaken}-\reff{switchshakenterminal}, see  Theorem 3.4 {\color{black} in} \cite{bmz} and the proof of  Theorem 2.3 of \cite{barlesjakobsen} for more details. More precisely {\color{black} by setting $\delta=\eps$},
 there exists a positive constant $C$ such that 
\b*
|v-v^{\eps,i}|_\infty\le C(\eps+k^\frac{1}{3}).
\e*
\no{\bf Step 2.} 
Let  $v^{(i)}_\eps:=v^{\eps,i}* \rho^\eps$, where $\{\rho^\eps\}$ is as in \reff{defoverueps}. 
As in Lemma 4.2 of \cite{bmz} and Lemma 3.4 of \cite{barlesjakobsen}  for $\eps\le \left(12\sup_i|v^{(i)}_\eps|_1\right)^{-1} k$, {\color{black} for  $i_0\in \text{argmin}_i v_\eps^{(i)}(t,x)$, the function $v^{(i_0)}_\eps$ is a supersolution to 
\be\label{supersolw}
-\Lc^X v_\eps^{(i_0)}(t,x)-F\left(t,x,v_\eps^{(i_0)},D v_\eps^{(i_0)},D^2v_\eps^{(i_0)}\right)\ge0.
\ee}
Moreover,  for any $(t,x)\in[0,T)\times\R^d$, 
we have $v^{(i_0)}_\eps(t,x)<v^{(i)}_\eps(t,x)+k$. 
Therefore, for all $i$ we have
\b*
(w_{\eps}-v)(t,x)= (v^{(i_0)}_{\eps}-v)(t,x)\le(v^{(i_0)}_{\eps}-v^{(i)}_{\eps})(t,x)+(v^{(i)}_{\eps}-v)(t,x)\le k+ C(\eps+k^\frac{1}{3}).
\e*
Choosing $k= C_1\eps$ with $C_1= 12\sup_i|v^{(i)}_\eps|_1$, one can write
\be\label{approxvweps}
(w_{\eps}-v)(t,x) \le C\eps^\frac{1}{3}.
\ee
\no{\bf Step 3.} {\color{black} By the definition of $w_\eps$, for any $(t,x)$ and $i_0\in \text{argmin}_i v_\eps^{(i)}(t,x)$ we have $w_\eps(t,x)=v_\eps^{(i_0)}(t,x)$ and $w_\eps\le v_\eps^{(i_0)}$ elsewhere. Therefore, $\TT_h[{w}_\eps](t,x)\le \TT_h[{v}_\eps^{(i_0)}]$. Moreover, since \reff{regunderweps} is  satisfied by $v_\eps^{(i_0)}$, by Lemma 3.22 of \cite{ftw}, one can conclude that 
\b*
 \Rc_h[{v}_\eps^{(i_0)}](t,x)
\!\!\!\! &:=&\!\!\!\!
 \frac{{v}_\eps^{(i_0)}(t,x)-\TT_h[{v}_\eps^{(i_0)}](t,x)}{h}
\! +\!
 \mathcal{L}^X{v}_\eps^{(i_0)}(t,x)+F(t,x,{v}_\eps^{(i_0)},D{v}_\eps^{(i_0)},D^2{v}_\eps^{(i_0)})\ge-Ch\eps^{-3}.
 \e*}
Therefore due to \reff{supersolw}, 
$\frac{{w}_\eps(t,x)-\TT_h[{w}_\eps](t,x)}{h}\ge -Ch\eps^{-3} $ holds true. 
By Proposition \ref{propmaxpri}, one can get
\be\label{approxvhweps}
(v^h-w_\eps)(t,x)\le Ch\eps^{-3}.
\ee
Now, \reff{approxvweps} and \reff{approxvhweps} yield
\b*
(v^h-v)(t,x)\le C(\eps^\frac{1}{3}+\eps^{-3}h).
\e*
 By minimizing on $\eps>0$, the desired lower bound is obtained. \ep

\begin{Remark}[Stochastic scheme]\label{montecarlo}{\rm
Scheme \reff{scheme} produces a deterministic approximate solution. However, in practice, we approximate the expectations in \reff{TT} based on a randomly generated set of sample paths of the process $\Xh$. As a result, the approximate solution is not  deterministic anymore. By  following the line of argument in Section 4 of \cite{ftw}, one can show the almost sure convergence of this stochastic approximate solution and even provide the same rate of convergence in $\L^p(\Omega,\P)$.

More precisely, assume that $\E$ is approximated by $\hat\E^N$ where $N$ denotes the number of sample paths. Suppose that for some $p\ge1$, there exist constants $C_b,\lambda,\nu>0$ such that $\left\Vert\hat\E^N[R]-\E[R]\right\Vert_p\le C_bh^{-\lb}N^{-\nu}$ for a suitable class of random variables $R$  bounded by $b$. By replacing $\E$ with $\hat\E^N$ in the scheme \reff{scheme}, one obtains a stochastic approximate solution $\hat v^h_N$. Then, if we choose $N=N_h$ which is chosen to satisfy $\lim_{h\to0}N_h^\nu h^{\lb+2}=\infty$, then under assumptions of Theorem \ref{thmconv} 
\b*
 \hat v^h_{N_h}(\cdot,\omega)\longrightarrow v
 &&
 \mbox{locally uniformly,}
\e*
for almost every $\omega$ 
where $v$ is the unique viscosity solution of \reff{equation}-\reff{terminal}. In addition, if  
\begin{equation}\label{eq:relbhandN}
 \lim_{h\to 0}N_h^\nu h^{\lb+\frac{21}{10}} > 0,
 \end{equation}
  we have that $ \Vert v-\hat v^h_{N_h}\Vert_p
 \le
 Ch^{1/10},$ 
under the assumptions of Theorem \ref{thmrateconv}. }
\end{Remark}

{\color{black}
\section{Numerical results}
\label{Sec:numres}
\subsection{Risk neutral pricing of geometric American put option}\label{subsec:lin}
We consider a geometric  American put option on  three risky assets each of which follows a Black-Scholes dynamics under risk neutral probability measure. The payoff of the option is given by $(K-\xi(T))_+$ where $K$ and $T$ are respectively the strike price and maturity and $\xi(t):=\prod_{i=1}^3S_i(t)$ where 
\b*
dS(t)
&=&
\text{diag}(S(t))(rdt+\text{diag}(\Sigma)\cdot dW(t).
\e*
Here $S(t)=(S_i(t))_{i=1}^3$ is the vector of asset prices, $W(t)=(W_i(t))_{i=1}^3$ is a 3-dimensional Brownian motion, $r$ is the risk free interest rate, and $\Sigma=(\sigma_i)_{i=1}^3$ where $\sigma_i$ is the volatility of the $i$th asset. \\
The price of this option at time $t$ and for asset price vector $s=(s_1,s_2,s_3)$ is given by 
\be\label{5dim}
v(t,s):=\sup\E\left[e^{-r(\tau-t)}(K-\xi(\tau))_+\Bigr|S(t)=s\right]
\ee
where the supremum is over all stopping times $\tau\in[t,T]$ adapted to the filtration generated by the 3-dimensional Brownian motion and  $\E$ is the risk neutral expectation. It is well-known  that function $v$ satisfies the following differential equation
\b*
0&=&
\min\left\{-\partial_tv-\frac12 \sum_{i=1}^3s_i^2\sigma_i^2\partial_{s_is_i}v-r\sum_{i=1}^3s_i\partial_{s_i}v+rv, v - g\right\}\\
v(T,s)
&=&
g(s).
\e*
where $g(s)=(K-\prod_{i=1}^3s_i)_+$.
We treat this linear equation as a fully nonlinear one by separating the linear second order operator into two parts. More precisely, for some $\sigma^2_0$, we choose the linear and nonlinear parts to be $\Lc^X\phi:=\frac{\sigma^2_0}{2}\sum_{i=1}^3s_i^2\sigma_i^2\partial_{s_is_i}\phi+r\sum_{i=1}^3s_i\partial_{s_i} \phi$ and $F(\cdot,\cdot,r\phi,D\phi,D^2\phi)=\frac{1-\sigma^2_0}{2}\sum_{i=1}^3s_i^2\sigma_i^2\partial_{s_is_i}\phi$, respectively. This leads to the choice of diffusion $X(t):=(X_i(t))_{i=1}^3$
\b*
dX(t)
&=&
{\sigma_0}\text{diag}(X(t))\text{diag}(\sigma)\cdot dW(t)
\e*
for the approximation scheme \reff{scheme}.
On the other hand the approximation of the second order derivatives in \reff{hermit} is given by 
\b*
x_i^2\partial_{x_ix_i}v(t,x)\approx\frac{1-\sigma_0^2}{2\sigma_0^2}\E[v(t+h,x+\sigma_0\text{diag}(x)\text{diag}(\sigma)\cdot W(h))\frac{W_i(h)^2-h}{h^2}],
\e*
where $x=(x_i)_{i=1}^3$. 
For the numerical implementation we choose the continuous-time interest rate $r=.03$, volatility of all assets $\sigma_i=.1$, $T=1$, $K=8$, and we evaluate the option at time $t=0$ at $s=(2,2,2)$.
{\small
\begin {table}
\caption {{\color{black} The simulation results for geometric American put option regarded as a nonlinear problem. \# of sample paths = 6 million, $N=$ \# of time steps, Time = time of the algorithm in seconds.}} \label{tbl:1}
{\color{black}\begin{center}
\begin{tabular}{|c|c|c|c|}
\hline
$N$     &  Time   & $\sigma_0^2=.9$ & $\sigma_0^2=1$ \\
\hline\hline
5&92&.301173&.326258\\
\hline
10&234&.309001&.334205\\
\hline 
15&360&.312642&.337974\\
\hline 
20&499&.314978&.340397\\
\hline 
40&1050&.320354&.347909\\
\hline
50&1325&.322115&.346041\\
\hline  
\end {tabular}
\end{center}}
\end {table}
}
The reference value for option price is obtained by applying the binomial tree algorithm to the 1-dimensional optimal stopping problem 
\b*
v(t,s):=\sup_\tau\E\left[e^{-r(\tau-t)}\xi(T)\Bigr|\xi(T)=\prod_{i=1}^3s_i\right].
\e*
 on the diffusion 
$\xi_t:=\prod_{i=1}^3S_i(t)$ satisfying
\b*
d\xi(t)
&=&
\xi(t)(3rdt+\bar\sigma dB_t)~~~~\text{with}~~~ ~\bar\sigma:=\Bigl(\sum_{i=1}^3\sigma_i^2\Bigr)^\frac12
\e*
where $B_t$ is a 1-dimensional Brownian motion. The binomial tree algorithm stabilizes to the value $0.338778$ for more than 20000 time steps.
\begin{figure}[h]
\centering{
\includegraphics[width=.7\textwidth]{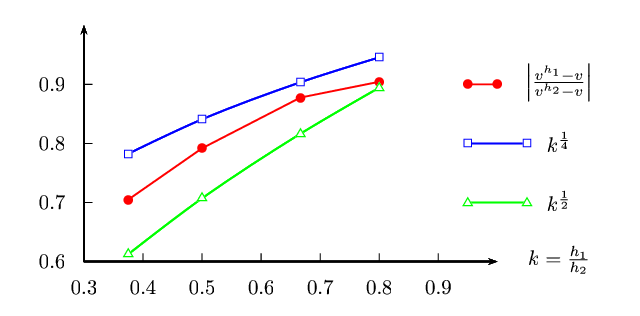}
  }
\caption{Rate of convergence analysis: $\bullet$: ratio of the error , i.e. $\left|\frac{v^{h_1}-v}{v^{h_2}-v}\right|$, $\vartriangle$: $\left({h_1}/{h_2}\right)^\frac{1}{2}$, and $\square$: $\left(h_1/{h_2}\right)^\frac{1}{4}$  (vertical axis) versus ${h_1}/{h_2}$ (horizontal axis).}\label{grph:1}
\end{figure}
The numerical result as well as the run time of the algorithm \footnote{\color{black} Dual core i5 2.5 GHz, 4 GB of RAM} is provided in Table \ref{tbl:1} for two different choices for $\sigma_0$. Here we used projection method described in \cite{bw} with $8^5$ locally linear basis functions with compact support. Notice that $\sigma^2_0=1$ corresponds to the Longstaff-Schwartz algorithm and for $\sigma^2_0=.9$, the inequality \reff{F2} is satisfied. In Figure \ref{grph:1} the red graph is the ratio of error  for two consecutive time steps plotted against the ratio of the time steps, while the green graph is the theoretical rate of convergence; i.e. $\frac{1}{4}$. Because of the analysis in \cite[Section 3.4]{ftw}, one expects to have a higher the rate of convergence for the scheme on the linear equations. Therefore, the simulated red plot must lie below the theoretical green plot. On the other hand, due to error of approximation  of expectations in the scheme, the rate of convergence $\frac12$ will never be obtained in practice.
\subsection{Indifference pricing of geometric American put option}
\label{subsec:nonlin}
We consider a geometric put option on two {\it non-tradable }risky assets with Black-Scholes dynamics given by
$S(t):=(S_1(t),S_1(t))$
\b*
dS_i(t)
&=&
S(t)(\mu_i dt+\sigma_i dW_i(t)),
\e*
where  $W(t)=(W_1(t),W_2(t))$ is a 2-dimensional Brownian motion, and $\mu_i$ and $\sigma_i$ are the drift  and   the volatility of the $i$th asset for $i=1, 2$ , respectively. We assume that there is a portfolio made of a tradable asset with Black-Scholes dynamics and money market with $r=0$ interest rate which satisfies
\b*
dX^\theta_t=\theta_t(\mu_0dt+\sigma_0dB(t)),
\e*
where $\theta$ is the amount of money in risky asset, $B(t)$ is a one dimensional Brownian motion, and $\mu_0$ and $\sigma_0$ are drift and volatility of the tradable asset, respectively. Here we assume that $dW_i(t)\cdot dB(t)=\rho_idt$ for $i=1, 2$.
The indifference pricing with exponential utility leads to the controller-stopper problem below
\begin{equation}\label{prob:3d}
\begin{split}
v(t,x,s_1,s_2):=\sup_{\tau,\theta}\E\left[-\exp(-\gamma(X^\theta(\tau)+(K-\prod_{i=1}^2S_i(\tau))_+))\Bigr|X(t)=x, S_i(t)=s_i,\; i=1, 2\right].
\end{split}
\end{equation}
which satisfies the fully nonlinear obstacle problem  below:
\b*
\min\left\{-\frac{\partial v}{\partial t}+\frac{(\mu v_x+\sum_{i=1}^2\sigma_0\rho_i\sigma_is_i\partial_{xs_i}v)^2}{2\sigma_0^2\partial_{xx}v}-\Lc^Sv,v-g\right\}
&=&0\\
v(T,x,s_1,s_2)= -\exp(-\gamma(x+(K-\prod_{i=1}^2s_i)_+))
\e*
where $\gamma>0$ is a constant and $g(t,x,s_1,s_2)=- \exp(-\frac{\mu_0^2}{2\sigma_0^2}(T-t)-\gamma(x+(K-\prod_{i=1}^2s_i)_+))$, $K$ is the strike price and $\Lc^S=\sum_{i=1}^2s_i\mu_i\partial_{s_i} +\frac12\sum_{i=1}^2s_i\sigma^2_i\partial_{s_is_i}$. 
To solve the above free boundary problem using scheme \reff{scheme}, we choose the linear and nonlinear parts as follow:
\b*
\bar\Lc\phi&=&\Lc^S\phi+\frac12\eps^2\partial_{xx}\phi\\
F(\cdot, \Dc\phi)&=&-\frac{(\mu \phi_x+\sum_{i=1}^2\sigma_0\rho_i\sigma_is_i\partial_{xs_i}\phi)^2}{2\sigma_0^2\partial_{xx}\phi}-\frac12\eps^2\partial_{xx}\phi.
\e*
\begin {table}
\caption {{\color{black} $M=$ \# of sample paths in million, $N=$ \# of time steps, Time = time of the algorithm on -dimensional problem in seconds, $\hat v$ and $\hat u$  value obtained by scheme \reff{scheme} for the 3-dimensional  problem \reff{prob:3d} and the 2-dimensional problem with the same value function. Notice that $v(0,1,1,1)=u(0,1,1)$.}} \label{tbl:nonlin}
{\color{black}\begin{center}
\begin{tabular}{|c|c|c|c|c|}
\hline
$N$ &  $M$        &  Time   & $\hat v(0,1,1,1)$ & $\hat u(0,1,1)$ \\
\hline\hline
5&2&34&-.341675&-.349489\\
\hline
\multirow{2}{*}{10}&1&36&-.303425&-.352110\\
\cline{2-5}
&2&76&-.356332&-.351678\\
\hline
\multirow{3}{*}{20}&2&180&-.356126&-.351659\\
\cline {2-5}
&3&273&-.351659&-.356126\\
\cline {2-5}
&4&372&-.348773&-.350201\\
\hline
30&3&422&-.353311&-.353088\\
\hline 
40&4&696&-.322095&-.361026\\
\hline 
\end {tabular}
\end{center}}
\end {table}
\begin{figure}[h]
\centering{
\includegraphics[width=.7\textwidth]{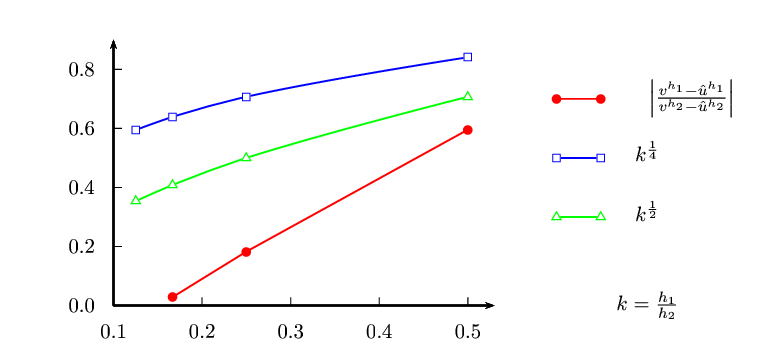}
  }
\caption{Rate of convergence analysis: $\bullet$: ratio of the error , i.e. $\left|\frac{v^{h_1}-\hat u^{h_1}}{v^{h_2}-\hat u^{h_2}}\right|$, $\vartriangle$: $\left({h_1}/{h_2}\right)^\frac{1}{2}$, and $\square$: $\left(h_1/{h_2}\right)^\frac{1}{4}$ (vertical axis) versus $k={h_1}/{h_2}$ (horizontal axis). } \label{grph:nonlin}
\end{figure}
Thus, the appropriate diffusion to be used inside \reff{scheme} is
\b*
d\bar X&=&\eps d\bar B(t),\\
dS(t)
&=&
\text{diag}(S(t))(\mu dt+\text{diag}(\Sigma)\cdot dW(t),
\e*
where $\mu=(\mu_1,\mu_2)$, $\Sigma=(\sigma_1,\sigma_2)$, and $\bar B(t)$ is a one dimensional Brownian motion independent of $W(t)$.\\
To find a reference value for the solution, we follow the same idea as in Section \ref{subsec:lin}. Since $\xi_t:=\prod_{i=1}^2S_i(t)$ satisfies
\b*
d\xi(t)
&=&
\xi(t)(\bar\mu dt+\bar\sigma dB_t)~~~~\text{with}~~~ ~\bar\mu:=\sum_{i=1}^2\mu_i,\;\bar\sigma:=\Bigl(\sum_{i=1}^2\sigma_i^2\Bigr)^\frac12,
\e*
we have $v(t,x,s_1,s_2)=u(t,x,\prod_{i=1}^2s_i)$ where the function $u$ is the solution of the 2-dimensional controller-stopper problem
\begin{equation}\label{prob:2d}
\begin{split}
u(t,x,s)
=
\sup_{\tau,\theta}\E\left[-\exp(-\gamma(X^\theta(\tau)+ (K-\xi(\tau))_+))\Bigr|X(t)=x, \xi(t)=s\right].
\end{split}
\end{equation}
Neither function $v$ nor $u$ have closed form solutions, but
we expect that if the scheme converges numerically, it approximates the function $u$ more accurately because of the reduction in the dimension. This is because the number of sample paths and time steps for a 2-dimensional problem for $u$ can be chosen larger than the 3-dimensional problem for $v$. Therefore, to examine the convergence of scheme, we compare the approximation of these functions by the scheme \reff{scheme}.
We set  $K=1$, $\gamma=1$, $T=1$ , $\eps=.05$, $X(0)=1$, $\rho_i=.1$, $\mu_0=\sigma_0=\mu_i=\sigma_i=.1$, and $S_i(0)=1$ for $i=1, 2$.\\
 The result of the simulation is summarized in Table \ref{tbl:nonlin}. In Figure \ref{grph:nonlin}, we establish the convergence analysis for the non-linear problem by using the approximations in Table \ref{tbl:nonlin} with the largest number of sample paths for each time step. 
}

\setcounter{section}{0}%
\renewcommand{\thesection}{\Alph{section}}%

\section{Appendix: A switching system with an obstacle}
\label{Sec:appendix}
In this section, we will provide some results needed in the Section \ref{proofrate}. In particular, we present a continuous dependence result for the switching system with obstacle and as a corollary a comparison result, which provides the uniqueness of the solution. Then, the existence and regularity of the solutions to the switching systems are provided. 

Consider the following system of PDEs for $v=(v^{(i)})_{i=1}^M$:
\begin{equation}\label{switchobstacle}
\begin{split}
\min\left\{\max\left\{-v^{(i)}_t-F_i(\cdot, v^{(i)}, Dv^{(i)},D^2 v^{(i)}),v^{(i)}-\Mc^{(i)}v\right\},v^{(i)}-g\right\}&=0,~\text{ for }i=1,\cdots,M;\\
v^{(i)}(T,\cdot)&=g(T,\cdot).
\end{split}
\end{equation}
We also need to consider  a variant of equation \reff{switchobstacle} as follows:
\begin{equation}\label{switchobstaclehat}
\begin{split}
\min\left\{\max\left\{-v^{(i)}_t-\hat F_i(\cdot, v^{(i)}, Dv^{(i)},D^2 v^{(i)}),v^{(i)}-\Mc^{(i)}v\right\},v^{(i)}-\hat g\right\}&=0,~\text{ for }i=1,\cdots,M;\\
v^{(i)}(T,\cdot)&=\hat g(T,\cdot).
\end{split}
\end{equation}

\no{\bf Assumption HJB-S.} {\it  We assume that in \reff{switchobstacle} and  \reff{switchobstaclehat}
\begin{align}\label{nonlin}
F_i(\cdot, v^{(i)}, Dv^{(i)},D^2 v^{(i)})=\inf_{\alpha\in\Ac^i}\Lc^{i,\alpha} v^{(i)} \; \text{ and } \;
\hat F_i(\cdot, v^{(i)}, Dv^{(i)},D^2 v^{(i)})=\inf_{\alpha\in\Ac^i}\hat\Lc^{i,\alpha} v^{(i)},
\end{align}
 $\Mc^{(i)}v=\min_{j:j\ne i}\{v_j+k\}$, $k$ is a non-negative constant, and
\begin{align*}
\Lc^{i,\alpha}\vp(x)&:=\frac12\Tr{a_i^\alpha(t,x)D^2\vp}+b_i^\alpha(t,x)D\vp+ c_i^\alpha(t,x)\vp+f_i^\alpha(t,x),\\
\hat\Lc^{i,\alpha}\vp(x)&:=\frac12\Tr{\ah_i^\alpha(t,x)D^2\vp}+\bh_i^\alpha(t,x)D\vp+ \ch_i^\alpha(t,x)\vp+\fh_i^\alpha(t,x). 
\end{align*}
Moreover,
\b*
L:= |g|_1+|\hat{g}|_1+\sup_{\alpha\in\bigcup_i\Ac^i}\left(|\sigma^\alpha|_1+|b^\alpha|_1+|c^\alpha|_1+|f^\alpha|_1+|\hat\sigma^\alpha|_1+|\bh^\alpha|_1+|\ch^\alpha|_1+|\fh^\alpha|_1\right) < \infty,
 \e*
and for all $\alpha\in\bigcup_i\Ac^i$, we have $c^\alpha_i, \ch^\alpha_i\le-1$.
}
\begin{Remark}
$c^\alpha_i, \ch^\alpha_i\le-1$ in Assumption {\bf HJB-S} is only to make the proofs simpler and is not a loss of generality.  This can be seen by applying the  change of variable $v^{(i)}\to e^{C(T-t)}v^{(i)}$ in the equations \reff{switchobstacle} and \reff{switchobstaclehat} for $C$ large enough.
\end{Remark}
\begin{Remark}\label{infsup}
For the sake of simplicity in Assumption {\bf HJB-S}, we only included the nonlinearities of infimum type. However, all the results of this appendix still hold if we assume that  
\be\label{P}
F_i(\cdot, v^{(i)}, Dv^{(i)},D^2 v^{(i)})=\inf_{\alpha\in\Ac^i}\sup_{\beta\in\Bc^i}\Lc^{i,\alpha,\beta}~~~~~~~~~\text{ and }~~~~~~~~~~
\hat F_i(\cdot, v^{(i)}, Dv^{(i)},D^2 v^{(i)})=\inf_{\alpha\in\Ac^i}\sup_{\beta\in\Bc^i}\hat\Lc^{i,\alpha,\beta},
\ee
\b*
 |g|_1+|\hat{g}|_1\!\!+\!\!\!\!\!\!\!\!\!\sup_{\tiny
 \begin{array}{c}
 \alpha\in\bigcup_i\Ac^i\\
 \beta\in\bigcup_i\Bc^i
 \end{array}}\!\!\!\!\!\left(|\sigma^{\alpha,\beta}_i|_1+|b^{\alpha,\beta}_i|_1+|c^{\alpha,\beta}_i|_1+|f^{\alpha,\beta}_i|_1+|\hat\sigma^{\alpha,\beta}_i|_1+|\bh^{\alpha,\beta}_i|_1+|\ch^{\alpha,\beta}_i|_1+|\fh^{\alpha,\beta}_i|_1\right) <\infty,
 \e*
and for all $\alpha\in\bigcup_i\Ac^i$ and  $\beta\in\bigcup_i\Bc^i$, we have  $c^{\alpha,\beta}_i, \ch^\alpha_i\le-1$. This remark is also valid if we change the order of $\inf$ and $\sup$ in \reff{P}.
\end{Remark}

\begin{Lemma}\label{noloop}
Let $u=(u^{(i)})_i$ and $v=(v^{(i)})_i$ be respectively the upper semicontinuous subsolution and the lower semicontinuous supersolution of \reff{switchobstacle} and  \reff{switchobstaclehat}, and  assume that $\vp(t,x,y)$ is a smooth function bounded from below. Define
\begin{align*}
\psi^{(i)}(t,x,y)&\equiv u^{(i)}(t,x)-v^{(i)}(t,y)-\vp(t,x,y),\\
\Jc_1&:=\left\{j\Bigl|~\exists (t',x',y'):~~\sup_{i,t,x,y}\psi^{(i)}(t,x,y)=\psi^{(j)}(t',x',y')\right\},\\
\Jc_2(t,x)&:=\left\{j\Bigl|~u^{(j)}(t,x)\le g(t,x)\right\}.
\end{align*}
Suppose that there exists an $(i'_0,t_0,x_0,y_0)$ such that  $\sup_{i,t,x,y}\psi^{(i)}(t,x,y)=\psi^{(i'_0)}(t_0,x_0,y_0)$and $\Jc_1\cap\Jc_2(t_0,x_0)=\emptyset$. Then, there exists an $i_0$ such that  $\psi^{(i_0)}(t_0,x_0,y_0)=\psi^{(i'_0)}(t_0,x_0,y_0)$ and 
\be\label{no_switch}
v^{(i_0)}(t_0,y_0)<\Mc^{(i_0)}v(t_0,y_0).
\ee
Moreover, if 
in a neighborhood of $(t_0,x_0,y_0)$ there are some continuous functions $h_0(t,x,y)>0$, $h(t,x)$ and $\hat h(t,y)$ such that
\b*
D^2\vp(t,x,y)\le h_0(t,x,y)\left(\begin{array}{cc}I&-I\\-I&I\end{array}\right)+\left(\begin{array}{cc}h(t,x)&0\\0&\hat h(t,y)\end{array}\right),
\e* 
then there are $a,b\in \R$ and $X,Y\in \S^d_+$ such that
\be\label{jets}
a-b&=&\vp_t(t_0,x_0,y_0),\\
\left(\begin{array}{cc}X&0\\0&-Y\end{array}\right)&\le& 2h_0(t_0,x_0,y_0)\left(\begin{array}{cc}I&-I\\-I&I\end{array}\right)+\left(\begin{array}{cc}h(t_0,x_0)&0\\0&\hat h(t_0,y_0)\end{array}\right),\label{jet2}
\ee
\begin{align}
-a-\inf_{\alpha\in\Ac_{i_0}}\left\{ \frac12\Tr{a_{i_0}^\alpha(t_0,x_0)X} +b_{i_0}^\alpha(t_0,x_0)D_x\vp(t_0,x_0,y_0)+ c_{i_0}^\alpha(t_0,x_0)u^{(i_0)}(t_0,x_0)+f_{i_0}^\alpha(t_0,x_0)\right\}&\le0,\label{sub}\\
-b-\inf_{\alpha\in\Ac_{i_0}}\left\{ \frac12\Tr{\ah_{i_0}^\alpha(t_0,y_0)Y} +\bh_{i_0}^\alpha(t_0,y_0)(-D_y\vp(t_0,x_0,y_0))+ \ch_{i_0}^\alpha(t_0,y_0)v^{(i_0)}(t_0,y_0)+\fh_{i_0}^\alpha(t_0,y_0)\right\}&\ge0.\label{super}
\end{align}
\end{Lemma}
\proof The first part of the proof is similar to those of Lemma A.2 of \cite{barlesjakobsen1}, Lemma A.1 of \cite{bmz}. The second part follows as a result of Theorem 2.2 of \cite{cil}.\ep\\
The following theorem on continuous dependence is used in Section \ref{subsec:thmrateconvii} and in the regularity result, Theorem \ref{switchreg} below.
Intuitively speaking, continuous dependence result asserts that a slight change in the coefficients of  \reff{switchobstacle}  changes the solution only slightly. 
\begin{Theorem}[Continuous dependence]\label{contdep}
Let {\bf HJB-S} hold.
Suppose that $u=(u^{(i)})_i$ and $v=(v^{(i)})_i$ are a bounded upper semicontinuous subsolution  of \reff{switchobstacle} and a bounded lower semicontinuous supersolution of  \reff{switchobstaclehat}, respectively. Then, for any $i=1,\cdots,N$, 
\begin{align*}
u^{(i)} -v^{(i)}
&\le\Bc:=C\max_j\Biggl\{
|(g-\hat g)(\cdot,\cdot)|_\infty+T\sup_\alpha\left\{|f^{j,\alpha}-\hat f^{j,\alpha}|_\infty+(|u|_\infty\vee|v|_\infty) |c^{j,\alpha}-\hat c^{j,\alpha}|_\infty\right\}\\
&+\sqrt{T}\sup_\alpha\left\{|\sigma^{j,\alpha}-\hat\sigma^{j,\alpha}|_\infty+|b^{j,\alpha}-\hat b^{j,\alpha}|_\infty\right\}\Biggr\}.
\end{align*}
\end{Theorem}
\proof
Let $\vp(t,x,y)=e^{\lb(T-t)}\frac{\theta}{2}|x-y|^2+e^{\lb(T-t)}\frac{\eps}{2}(|x|^{2}+|y|^{2})$
and define 
$$\Dc := \sup_{t,i,x,y}\left\{u^{(i)}(t,x)-v^{(i)}(t,y)-\vp(t,x,y)-\frac{\bar \eps }{t}\right\},$$
 where  $\eps, \bar\eps>0$ are arbitrary constants and constants $\lb, , \theta>0$ will be determined later in the proof.
We will show that $\Dc$ is bounded by a constant $B(\eps,\bar\eps,\theta)$ which is bounded by $\Bc$ mentioned in the theorem as {\color{black}$(\eps,\bar\eps)\to(0,0)$ and $\theta$ is set appropriately.} Then, it would follow that
\b*
u^{(i)}(t,x)-v^{(i)}(t,x)\le \Dc+\frac{\bar\eps}{t}+e^{\lb(T-t)}\frac{2\eps|x|^{2}}{2}\le B(\eps,\bar\eps,\theta)+\frac{\bar\eps}{t}+e^{\lb(T-t)}\frac{2\eps|x|^{2}}{2}.
\e*
Sending $\eps, \bar\eps\to 0$, one would then obtain
\b*
u^{(i)}(t,x)- v^{(i)}(t,x)\le\Bc,\text{   for $t>0$}.
\e*
Note that the above inequality is also valid for $t=0$ by considering $[-\delta,T]$ as the time interval {\color{black} and} by changing $T$ to $T+\delta$.\\
  Define 
$$\psi^{(i)}(t,x,y)=u^{(i)}(t,x)-v^{(i)}(t,y)-\vp(t,x,y)-\frac{\sigma(T-t)}{2T}-\frac{\bar \eps }{t},$$
where $\sigma=\Dc-\sigma_T$ with
$\sigma_T=\sup_{i,x,y}\left\{u^{(i)}(T,x)-v^{(i)}(T,y)-\vp(T,x,y)-\frac{\bar \eps }{T}\right\}^+$.
Let
\be\label{suppenalty}
\bar\Dc:=\sup_{i,t,x,y}\psi^{(i)}(t,x,y).
\ee
Since $u^{(i)}$ and $v^{(i)}$ are bounded, we have $\bar\Dc<\infty$. On the other hand, by semicontinuity of $u^{(i)}$ and $v^{(i)}$, one can conclude that the supremum in the definition of $\bar\Dc$ is attained at some point $(i_0,t_0,x_0,y_0)$. In other words, $\Jc_1\neq\emptyset$ (see Lemma \ref{noloop} for the definition of $\Jc_1$).

If $\sigma\le0$, then $\Dc\le \sigma_T.$
Since
\b*
\sigma_T&\le&|g-\hat g|_\infty+\sup_{x,y}\{|g|_1|x-y|-\vp(T,x,y)\}-\frac{\bar\eps}{T}\\
&\le&|g-\hat g|_\infty+
\sup_{x,y}\{|g|_1|x-y|-\frac{\theta}{2}|x-y|^2\}\le|g-\hat g|_\infty+
\frac{|g|_1^2}{2\theta},
\e*
one can conclude that $\Dc\le |g-\hat g|_\infty+\frac{|g|_1^2}{2\theta}$. Therefore, we may assume that $\sigma>0$.
From the definition of $\bar\Dc$, we  have $t_0>0$. On the other hand,  $\sigma>0$ implies $t_0<T$. Because if $t_0=T$, then $\sigma_T\ge\bar\Dc$ which implies
\b*
\sigma_T
&\ge&
\Dc-\frac{\sigma}{2}\ge\sigma_T+\frac{\sigma}{2}>\sigma_T
\e*
which is a contradiction.  So, we have $0<t_0<T$. We continue the proof by considering two different cases.

\no{\bf Case 1: $\Jc_1\cap\Jc_2(t_0,x_0)\neq\emptyset$}. The supremum in \reff{suppenalty} is attained at some point $(i_0,t_0,x_0,y_0)$ with
$u^{(i)}(t_0,x_0)\le g(t_0,x_0)$ and $v^{(i)}(t_0,y_0)\ge \hat g(t_0,y_0)$.
Therefore, 
\b*
\bar\Dc
&\le& 
g(t_0,x_0)-\hat g(t_0,y_0)-\vp(t_0,x_0,y_0)-\frac{\sigma \;(T-t_0)}{2T}-\frac{\bar \eps }{t_0}\\
&\le&
|g-\hat g|_\infty+|g|_1|x_0-y_0|-\frac{\theta}{2}|x_0-y_0|^2\le |g-\hat g|_\infty+\frac{|g|^2_1}{2\theta}.
\e*
On the other hand, since 
$\Dc\le\bar\Dc+\frac{\sigma}{2} \le \bar\Dc+ \frac12(\Dc-\sigma_T)\le \bar\Dc + \frac12\Dc$, 
we have $\Dc\le 2 |g-\hat g|_\infty+ \frac{|g|^2_1}{\theta}$.

\no{\bf Case 2: $\Jc_1\cap\Jc_2(t_0,x_0)=\emptyset$}. 
In this case, {\color{black} \reff{no_switch} is satisfied and by Lemma \ref{noloop} and   the same line of argument as Theorem A.1 in  \cite{jakobsen} the result is provided. For the convenience of the reader, we present a sketch of the proof.

By subtracting \reff{super} from \reff{sub}, we have 
\b*
\lb\phi(t_0,x_0,y_0)+\frac{\Dc-\sigma_T}{2T}+\frac{\bar\eps}{t_0^2} 
&-&
\inf_{\alpha\in\Ac_{i_0}}\Bigl\{ \frac12(\Tr{a_{i_0}^\alpha(t_0,x_0)X}-\Tr{\ah_{i_0}^\alpha(t_0,y_0)Y})\\ 
&-&
b_{i_0}^\alpha(t_0,x_0)D_x\vp(t_0,x_0,y_0)- \bh_{i_0}^\alpha(t_0,y_0)D_y\vp(t_0,x_0,y_0) \\
&-&
c_{i_0}^\alpha(t_0,x_0)(u^{(i_0)}(t_0,x_0)-v^{(i_0)}(t_0,y_0))\\
&-& 
(c_{i_0}^\alpha(t_0,x_0)- \ch_{i_0}^\alpha(t_0,y_0))v^{(i_0)}(t_0,y_0)\\
&-&
f_{i_0}^\alpha(t_0,x_0)-\fh_{i_0}^\alpha(t_0,y_0)\Bigr\}\le 0
\e*
Now, using $ c_{i_0}^\alpha\le -1$ and $u^{(i_0)}(t_0,x_0)-v^{(i_0)}(t_0,y_0)\ge \bar \Dc\ge \frac12 \Dc>0$ together with \reff{jets}-\reff{jet2}, one can obtain the following bound for $\Dc$
\b*
\Dc
&\le &
CT\Bigl(\theta \sup_\alpha\left\{|\sigma^{j,\alpha}-\hat\sigma^{j,\alpha}|^2_\infty+|b^{j,\alpha}-\hat b^{j,\alpha}|^2_\infty\right\}\\
&&+\sup_\alpha\left\{|f^{j,\alpha}-\hat f^{j,\alpha}|_\infty+(|u|_\infty\vee|v|_\infty) |c^{j,\alpha}-\hat c^{j,\alpha}|_\infty\right\}+\sigma_T
\Bigr)+C_1|x_0-y_0|^2-\lb\phi(t_0,x_0,y_0),
\e*
where $C_1$ is the constant depending only on $L$ in Assumption {\bf HJB-S}. After we choose $\lb\ge C_1$ in the above and maximize the right hand side with respect to $\theta$, the proof is complete.
}\ep

The following result {\color{black} is} a straightforward consequence of Theorem \ref{contdep}, and will be used to establish the existence and the regularity of the solution to \reff{switchobstacle}.
\begin{Corollary}\label{comparison}
Assume that {\bf HJB-S} holds.
Suppose that $u=(u^{(i)})_i$ and $v=(v^{(i)})_i$ are respectively a bounded upper semicontinuous subsolution and a bounded lower semicontinuous supersolution of  \reff{switchobstacle}. Then, for any $i=1,\cdots,N$, $u^{(i)}\le v^{(i)}$.
\end{Corollary}
\begin{Theorem}[Existence]\label{switchexist}
Assume that {\bf HJB-S} holds. Then there exists a unique continuous viscosity solution in the class of bounded functions to \reff{switchobstacle}.
\end{Theorem}
\proof We follow Perron's method (see e.g. Section 4 of \cite{cil}). Observe that by Assumption {\bf HJB-S}, $\underline{u}=-K$ and $\overline{v}=K$ are respectively sub and supersolutions of  \reff{switchobstacle} for a suitable choice of  positive constant $K$. Define
$
v^{(i)}(t,x):=\sup\{u^{(i)}(t,x)\;;\; u\text{ is a subsolution to \reff{switchobstacle} }\}
$
and
\b*
v^{(i)*}(t,x):={\color{black}\lim_{\delta \to 0}} \sup\{v^{\color{black}(i)}(s,y): |x-y|+|s-t| \leq \delta, s \in [0,T] \},
\e*
and
\b*
v^{(i)}_*(t,x):={\color{black}\lim_{\delta \to 0}} \inf\{v^{\color{black}(i)}(s,y): |x-y|+|s-t| \leq \delta, s \in [0,T] \}.
\e*
It is straight forward that  $-K\le v^{(i)}_*\le v^{(i)*}\le K$.
We want to show that $(v^{(i)*})_{i=1}^M$ and $(v^{(i)}_*)_{i=1}^M$ are respectively a sub and a supersolution to \reff{switchobstacle}  which by comparison, Corollary \ref{comparison}, yields the desired result.

\no{\bf Step 1: Subsolution property of $v^{(i)*}$.}
We start by showing that $(U,\cdots,U)$ with
\b*
U(t,x):=a_\eps(T-t)+g(T,z)+|g|_1\left(T-t+|x-z|^2+\eps\right)^\frac12
\e*
is a supersolution to\reff{switchobstacle}  for a suitable positive constant $a_\eps$. Observe that since
\b*
U(t,x)-g(t,x)\ge g(T,z) -g(t,x)+|g|_1\left(T-t+|x-z|^2+\eps\right)^\frac12\ge0,
\e*
we have that $U(t,x)\ge g(t,x)$,  and in particular $U(T,x)\ge g(T,x)$.  On the other hand, by simple calculations, one can show that, for an appropriate choice of $a_\eps$, we have $-U_t-\inf_{\alpha\in\Ac^i}\Lc^{i,\alpha} U\ge0$.

Therefore, by comparison, Corollary  \ref{comparison}, for any subsolution $u$, $u\le U$ which implies $v^{(i)*}\le U$; specially $v^{(i)*}(T,x)\le U(T,x)$. Sending $\eps\to0$ and setting $x=z$, $v^{(i)*}(T,x)\le g(T,x)$.

Now, for fixed $i$, we suppose $t<T$ and $\vp$ is a test function such that
\b*
0=\max_{[0,T]\times\R^d}\{v^{(i)*}-\vp\}=(v^{(i)*}-\vp)(t,x).
\e*
It follows from the definition of $v^{(i)*}$ that there exists a sequence $\{(u_n,t_n,x_n)\}_n$ with  $t_n<T$  such that $u_n$ is a subsolution to \reff{switchobstacle}, $(t_n,x_n)\to(t,x)$, $u^{(i)}_n(t_n,x_n)\to v^{(i)*}(t,x)$, and $(t_n,x_n)$ is the global strict maximum of $u^{(i)}_n-\vp$. 
Let $\delta_n:=\max_{[0,T]\times\R^d}\{u^{(i)}_n-\vp\}$. By the subsolution property of $u_n$, we have
\b*
\min\left\{\max\left\{-\vp_t-\inf_{\alpha\in\Ac^i}\Lc^{i,\alpha} (\vp+\delta_n),u_n^{(i)}-\Mc^{(i)}u_n\right\},\vp+\delta_n-g\right\}(t_n,x_n)\le0.
\e*
Because $\Mc^{(i)}u_n\le\Mc^{(i)}v^*$, by sending $n\to\infty$,
\b*
\min\left\{\max\left\{-\vp_t-\inf_{\alpha\in\Ac^i}\Lc^{i,\alpha} \vp,v^{(i)*}-\Mc^{(i)}v^*\right\},v^{(i)*}-g\right\}(t,x)\le0.
\e*

\no{\bf Step 2: Supersolution property of $v^{(i)}_*$.}
Since $(g,\cdots,g)$ is a subsolution to \reff{switchobstacle}, $v^{(i)}_*(t,x)\ge g(t,x)$. In particular, $v^{(i)}_*(T,x)\ge g(T,x)$. Therefore, we only need to show that 
\be\label{superv*}
\max\left\{-(v^{(i)}_*)_t-\inf_{\alpha\in\Ac^i}\Lc^{i,\alpha} v^{(i)}_*,v^{(i)}_*-\Mc^{(i)}v_*\right\}\ge0,
\ee
 on $[0,T)\times\R^d$ in the viscosity sense. We will prove \eqref{superv*} by a contradiction argument. Assume that there are a test function $\vp$ and $(i,t,x)$ with $t<T$ such that $(t,x)$ is the global strict minimum of $v^{(i)}_*-\vp$ and {\color{black} $(v^{(i)}_*-\vp)(t,x)=0$ but} $\max\left\{-\vp_t-\inf_{\alpha\in\Ac^i}\Lc^{i,\alpha}\vp ,\vp-\Mc^{(i)}v_*\right\}(t,x)<0$.
Then, by continuity of $\vp$ and the equation and lower semicontinuity of $v_*$, one can find $\eps>0$ and $\delta>0$ small enough, such that for $|x-y|+|s-t|<\delta$ we have that $\vp+\eps<v_*^{(i)}$ and that
\be\label{vpeps}
\max\left\{-(\vp+\eps)_t-\inf_{\alpha\in\Ac^i}\Lc^{i,\alpha}(\vp+\eps) ,(\vp+\eps)-\Mc^{(i)}v_*\right\}(s,y)<0.
\ee
Define 
\b*
w^{(j)}(s,y):=\left\{\begin{array}{ll}\max\{\vp+\eps,v^{(j)*}\}(s,y), &j=i ~ \text{ and } ~ |x-y|+|s-t|<\delta;\\v^{(j)*}(s,y),&\text{otherwise.}\end{array}\right.
\e*
Since $v^*$ is a subsolution to \reff{switchobstacle} and by \reff{vpeps}, one can show that $w$ is a subsolution to \reff{switchobstacle}. By the definition of $v_*^{(i)}$, we must have $v_*^{(i)}\ge w^{(i)}$, which contradicts the fact that $w^{(i)}(t,x)=\vp(t,x)+\eps<v_*^{(i)}(t,x)$ for $|x-y|+|s-t|<\delta$. \ep
\begin{Theorem}[Regularity]\label{switchreg}
Assume that {\bf HJB-S} holds.
Let $(u^{(i)})_{i=1}^M$ be the solution to \reff{switchobstacle}. Then, $(u^{(i)})_{i=1}^M$ is Lipschitz continuous with respect to $x$ and $\frac12$-H\"older continuous with respect to $t$ on $\R^d\times[0,T]$.
\end{Theorem}
\proof
{\bf Lipschitz continuity with respect to $x$:}  For fixed $y\in\R^d$, $v^{(i)}(x)=u^{(i)}(t,x+y)$ is the solution of a switching system obtained from \reff{switchobstacle} by replacing $\Lc^{i,\alpha}$ with
\b*
\Lc^{i,\alpha,y}\vp(x)&:=&\frac12\Tr{a_i^\alpha(t,x+y)D^2\vp}+b_i^\alpha(t,x+y)D\vp+ c_i^\alpha(t,x+y)\vp+f_i^\alpha(t,x+y),
\e* 
with the terminal condition given by $v^{(i)}(T,x)=g(T,x+y)$.
By Theorem \ref{contdep}, there is a positive constant $C$ such that
\b*
\sup_{t,x}|u^{(i)}(t,x)-u^{(i)}(t,x+y)|=\sup_{t,x}|u^{(i)}(t,x)-v^{(i)}(t,x)|&\le& C|y|.
\e*
{\bf $\frac12$-H\"older continuity with respect to $t$:} For $t<s$, define $\bar u=(\bar u^{(i)})_{i=1}^M$ to be the solution to
\begin{equation*}
\begin{split}
\max\left\{-\bar u^{(i)}_t-F_i(\cdot, \bar u^{(i)}, D\bar u^{(i)},D^2 \bar u^{(i)}),\bar u^{(i)}-\Mc^{(i)}\bar u\right\}&=0,~\text{ for }i=1,\cdots,M;\\
\bar u^{(i)}(s,\cdot)&=u^{(i)}(s,\cdot).
\end{split}
\end{equation*}
Since $\bar u$ is a subsolution of \reff{switchobstacle} on $[0,s]\times\R^d$ with terminal condition $u^{(i)}(s,\cdot)$, by comparison result, Corollary \ref{comparison}, we have $\bar u^{(i)}\le u^{(i)}$. Therefore, $u^{(i)}(t,x)-u^{(i)}(s,x)\ge\bar u^{(i)}(t,x)-\bar u^{(i)}(s,x)$. By Theorem A.1 of \cite{barlesjakobsen}, $\bar u^{(i)}$ is $\frac12$-H\"older continuous in $t$ which provides
\b*
 u^{(i)}(t,x)- u^{(i)}(s,x)\ge-C\sqrt{s-t}.
\e*
Now, for fixed $y\in\R^d$, define 
\b*
\psi^{(i)}(t,x):=\frac{\lambda L}{2}e^{A(s-t)}\left(|x-y|^2+B(s-t)\right)+\frac{L}{\lambda}+B(s-t)+g(s,y),
\e*
where $A$,  $B$ and $\lb$ are positive constants which will be given later and $L$ is the same as in Assumption {\bf HJB-S}. We will show that for an appropriate choice of $A$ and $B$, $(\psi^{(i)})_{i=1}^M$ is a supersolution of \reff{switchobstacle} with terminal condition $g(s,x)$. Then, comparison, Corollary~\ref{comparison}, would then imply that $u^{(i)}\le\psi^{(i)}$. Therefore,
\b*
\nonumber u^{(i)}(t,y)-u^{(i)}(s,y)&\le& \psi^{(i)}(t,y)-g(s,y)\le\frac{\lambda L}{2}e^{A(s-t)}B(s-t)+\frac{L}{\lambda}+B(s-t).
\e*
By setting $\lb=\frac{1}{\sqrt{s-t}}$, we have $
u^{(i)}(t,y)-u^{(i)}(s,y) \le C\sqrt{s-t}$,
where $C$ is a positive constant.

Therefore, it remains to show that for $A$ and $B$ large enough, we have 
\b*
\min\left\{\max\left\{-\psi^{(i)}_t-\inf_{\alpha\in\Ac^i}\Lc^{i,\alpha}\psi^{(i)} ,\psi^{(i)}-\Mc^{(i)}\psi^{(i)}\right\},\psi^{(i)}-g\right\}\ge0,
\e*
 on $[0,s]\times\R^d$.
Since $\psi^{(i)}-\Mc^{(i)}\psi^{(i)}<0$, one needs to show that 
\b*
-\psi^{(i)}_t-\inf_{\alpha\in\Ac^i}\Lc^{i,\alpha}\psi^{(i)} \ge0&\text{and}&\psi^{(i)}-g\ge0.
\e*
Observe that if $B\ge1$, by the regularity assumption on $g$, we have
\b*
\psi^{(i)}(t,x)-g(t,x)&\ge&\frac{L}{2}\left(\lb|x-y|^2+\lb(s-t)+\frac{2}{\lb}\right)+g(s,y)-g(t,x)\ge0.
\e*
On the other hand,
\begin{align*}
-\psi^{(i)}_t&-\inf_{\alpha\in\Ac^i}\Lc^{i,\alpha}\psi^{(i)} =
\sup_{\alpha\in\Ac^i}\Biggl\{
\frac{L\lb}{2}e^{A(s-t)}\biggl(A|x-y|^2+AB(s-t)+B\\
&-\frac12\Tr{a^{\alpha,i}}-b^{\alpha,i}\cdot(x-y)\biggr)+B-c^{\alpha,i}\psi^{(i)}-f^{\alpha,i}
\Biggr\}\\
&\ge\frac{L\lb}{2}e^{A(s-t)}\left(A|x-y|^2-L|x-y|+LB-\frac{L}{2}\right)+B-CL.
\end{align*}
By choosing $A$ and $B$ large enough, the right hand side in the above inequality is positive which completes the argument.\ep

\section{Appendix: Proof of  Lemma \ref{t-to-T}}
\label{Sec:appendix_b}
If  $v^h(t_i,x)=g(t_i,x)$ holds true, then since function $g$ is $\frac12$-H\"older continuous on $t$, the proof is done. So, we assume that 
$v^h(t_i,x)>g(t_i,x)$.  We introduce the discrete stopping time $\tauh:=\min\{t_j|j\ge i,~v^h(t_j,\Xh_{t_j})=g(t_j,\Xh_{t_j})\}$. Observe that $t_i<\tauh$.\\
{\bf Step 1.} Let $\bar\theta$ be such that $F(t_j,\Xh_{t_j}^x,\Dc_h v^h(t_j,X_{t_j}^x))-F(t_j,\Xh_{t_j}^x,0,0,0)=\nb F(t_j,\Xh_{t_j}^x,\bar\theta)\cdot\Dc_h v^h(t_j,X_{t_j}^x))$
For all $j=i,\cdots,n-1$, on the event $\{t_j<\tauh\}$ one can write
\b*
v^h(t_j,x)&=&\E_{t_j,x}[v^h(t_{j+1},\Xh_{t_{j+1}}^x)P_{j+1}]+hF_j,
\e*
where $F_j=F(t_j,\Xh_{t_j}^x,0,0,0)$, $\Delta W_{j+1}=W_{t_{j+1}}-W_{t_{j}}$, $P_{j+1}=1-\alpha_j+\sqrt{h}\beta_j\cdot\Delta W_{j+1}+h^{-1}\alpha_j\cdot\Delta W_{j+1}\Delta W_{j+1}^\text{T}$, and $\alpha_j:=F_\gamma\cdot a^{-1}(t_j,\Xh_{t_j}^x,\bar\theta)<1$ and $\beta_j:=F_p\cdot \sigma^{-1}(t_j,\Xh_{t_j}^x,\bar\theta)$ are $\Fc_{t_{j+1}}$-measurable. We can rewrite the above equality in the following form.
\begin{align}\label{1-step-ahead}
v^h(t_j,x)\1_{\{t_j<\tauh\}}&=&\E_{t_j,x}[v^h(t_{j+1},\Xh_{t_{j+1}}^x)P_{j+1}\1_{\{t_{j+1}<\tauh\}}+g(\tauh,\Xh_\tauh^x)P_{j+1}\1_{\{t_{j+1}=\tauh\}}]+h\1_{\{t_j<\tauh\}}F_j
\end{align}
Notice that the first term in the  right hand side of  \reff{1-step-ahead} is zero if $\{t_j\ge \tauh\}$. 
We define $Q_j:=\prod_{k=i+1}^jP_k$ with $Q_i:=1$ and $V_j:=v^h(t_{j},\Xh_{t_{j}}^x)Q_{j}\1_{\{t_{j}<\tauh\}}$. Observe that $Q_j$ is a discrete martingale with respect to $\{W_{t_j}\}_{j=i}^n$. Multiplying \reff{1-step-ahead} by $Q_j$, one can write
\b*
V_j=\E_{t_j,x}[V_{j+1}+g(\tauh,\Xh_\tauh^x)Q_{j+1}\1_{\{t_{j+1}=\tauh\}}]+h\1_{\{t_j<\tauh\}}Q_jF_j.
\e*
By summing the above equality over $j=i,\cdots, n-1$ and taking expectation $\E_{t_i,x}$, we have
\b*
v^h(t_i,x)&=&V_i=\E_{t_i,x}\left[V_{n}+\sum_{j=i}^{n-1}g(\tauh,\Xh_\tauh^x)Q_{j+1}\1_{\{t_{j+1}=\tauh\}} +h\sum_{j=i}^{n-1}\1_{\{t_j<\tauh\}}Q_jF_j\right]\\
&=&\E_{t_i,x}\left[\sum_{j=i}^{n-1}g(\tauh,\Xh_\tauh^x)Q_{j+1}\1_{\{t_{j+1}=\tauh\}} +h\sum_{j=i}^{n-1}\1_{\{t_j<\tauh\}}Q_jF_j\right].
\e*
Observe that here we used $V_n=0$ by the definition of $\tauh$. Thus, we can write
\be\label{vh-g}
v^h(t_i,x)-g(t_i,x)=\E_{t_i,x}\left[(g(\tauh,\Xh_\tauh^x)-g(t_i,x))\sum_{j=i}^{n-1}Q_{j+1}\1_{\{t_{j+1} =\tauh\}}+h\sum_{j=i}^{n-1}\1_{\{t_j<\tauh\}}Q_jF_j\right],
\ee
where in the above we used optional stopping theorem for  $\E_{t_i,x}[\sum_{j=i}^{n-1}Q_{j+1}\1_{\{t_{j+1} =\tauh\}}]=1$.
Our  goal is to show that the right hand side of \reff{vh-g} is bounded by $C\sqrt{T-t_i}$. First observe that by Assumption {\bf F}(i), $F_j$ is bounded. Then, because $Q_j$ is a positive martingale, we can bounded the second term in \reff{vh-g} by $C(T-t_i)$:
\b*
\left|\E_{t_i,x}\left[h\sum_{j=i}^{n-1}\1_{\{t_j<\tauh\}}Q_jF_j\right]\right|\le  C h\sum_{j=i}^{n-1}\E_{t_i,x}[Q_j]\le C(T-t_i).
\e*
We continue by bounding the other term in \reff{vh-g} in the next step.

\no{\bf Step 2.} To bound $\E_{t_i,x}\left[(g(\tauh,\Xh_\tauh^x)-g(t_i,x))\sum_{j=i}^{n-1}Q_{j+1}\1_{\{t_{j+1} =\tauh\}}\right]$, we want to apply It\^o formula on $g(\tauh,\Xh_\tauh^x)-g(t_i,x)$. But, because $g$ is not a $C^2$ function, we first approximate $g$ by a smooth function uniformly, i.e $|g-g_\eps|_\infty\le C\eps$. This can be done by $g_\eps:=g*\rho_\eps$  where $\{\rho_\eps\}_\eps$ is a family of mollifiers. Because $g$ is Lipschitz on $x$ and $\frac12$-H\"older on $t$, we have
\be\label{g-eps-bnd}
|\partial_tg_\eps|_\infty\le  {\eps}^{-1},~~~|Dg_\eps|_\infty\le C,~~~\text{and}~~~|D^2g_\eps|_\infty\le {\eps}^{-1}.
\ee
Therefore, we write
\b*
g(\tauh,\Xh_\tauh^x)-g(t_i,x)
&=&
(g_\eps(\tauh,\Xh_\tauh^x)-g_\eps(t_i,x))+(g(\tauh,\Xh_\tauh^x)-g_\eps(\tauh,\Xh_\tauh^x)) + (g_\eps(t_i,x)-g(t_i,x)).
\e*
Observe that since $|g_\eps-g|_\infty\le C\eps$, one has
\begin{equation}\label{g-eps-bnd}
\begin{split}
\left|\E_{t_i,x}\left[(g_\eps(\tauh,\Xh_\tauh^x)-g(\tauh,\Xh_\tauh^x))\sum_{j=i}^{n-1}Q_{j+1}\1_{\{t_{j+1} =\tauh\}}\right]\right|
&\le
C\eps,\\
\left|\E_{t_i,x}\left[(g_\eps(t_i,x)-g(t_i,x))\sum_{j=i}^{n-1}Q_{j+1}\1_{\{t_{j+1} =\tauh\}}\right]\right|
&\le
C\eps.
\end{split}
\end{equation}
In the following steps, we find a bound on $\E_{t_i,x}\left[(g_\eps(\tauh,\Xh_\tauh^x)-g_\eps(t_i,x))\sum_{j=i}^{n-1}Q_{j+1}\1_{\{t_{j+1} =\tauh\}}\right]$ in terms of $T-t_i$ and $\eps$.\\
{\bf Step 3.}  We apply It\^o formula on $g_\eps(\tauh,\Xh_\tauh^x)-g_\eps(t_i,x)$:
\b*
g_\eps(\tauh,\Xh_\tauh^x)-g_\eps(t_i,x)=\int_{t_i}^\tauh \Lc^\Xh g_\eps(s,\Xh_s^x)ds+\int_{t_i}^\tauh Dg_\eps(s,\Xh_s^x)\cdot dW_s,
\e*
where $\Lc^\Xh$ is the infinitesimal generator for the processe $\Xh$. Thus,
\begin{equation}\label{summation-until-tau}
\begin{split}
\E_{t_i,x}&\left[(g _\eps(\tauh,\Xh_\tauh^x)-g_\eps(t_i,x))\sum_{j=i}^{n-1}Q_{j+1}\1_{\{t_{j+1} =\tauh\}}\right]\\
&=
\sum_{j=i}^{n-1}\E_{t_i,x}\left[\left(\int_{t_i}^\tauh \Lc^\Xh g_\eps (s,\Xh_s^x)ds+\int_{t_i}^\tauh Dg_\eps (s,\Xh_s^x)\cdot dW_s\right)Q_{j+1}\1_{\{t_{j+1} =\tauh\}}\right].
\end{split}
\end{equation}
We proceed by calculating  the term in the above summation for $j=n-1$. 
\begin{equation}\label{bad-ass-term}
\begin{split}
\E_{t_i,x}&\left[\left(\int_{t_i}^\tauh \Lc^\Xh g_\eps (s,\Xh_s^x)ds+\int_{t_i}^\tauh Dg_\eps (s,\Xh_s^x)\cdot dW_s\right)Q_{n}\1_{\{t_{n} =\tauh\}}\right]\\
&=
\E_{t_i,x}\left[\left(\int_{t_i}^{t_n} \Lc^\Xh g_\eps (s,\Xh_s^x)ds+\int_{t_i}^{t_n} Dg_\eps (s,\Xh_s^x)\cdot dW_s\right)Q_{n}\1_{\{t_{n} =\tauh\}}\right]\\
&=
\E_{t_i,x}\left[\left(\int_{t_i}^{t_{n-1}} \Lc^\Xh g_\eps (s,\Xh_s^x)ds+\int_{t_i}^{t_{n-1}} Dg_\eps (s,\Xh_s^x)\cdot dW_s\right)Q_{n}\1_{\{t_{n} =\tauh\}}\right]\\
&+\E_{t_i,x}\left[\left(\int_{t_{n-1}}^{t_n} \Lc^\Xh g_\eps (s,\Xh_s^x)ds+\int_{t_{n-1}}^{t_n} Dg_\eps (s,\Xh_s^x)\cdot dW_s\right)Q_{n}\1_{\{t_{n} =\tauh\}}\right]
\end{split}
\end{equation}
We first bound the second term in the right hand side in the next step.\\
\no{\bf Step 4.} Since  $Q_{t_{n-1}}$ and $\1_{\{t_{n} =\tauh\}}$ are $\Fc_{t_{n-1}}$ measurable, the second term in the right hand side can be written as
\b*
\E_{t_i,x}\left[\left(\int_{t_{n-1}}^{t_n} \Lc^\Xh g_\eps (s,\Xh_s^x)ds+\int_{t_{n-1}}^{t_n} Dg_\eps (s,\Xh_s^x)\cdot dW_s\right)Q_{n}\1_{\{t_{n} =\tauh\}}\right]\\
=
\E_{t_i,x}\left[Q_{n-1}\1_{\{t_{n} =\tauh\}}\E_{t_{n-1}}\left[\left(\int_{t_{n-1}}^{t_n} \Lc^\Xh g_\eps (s,\Xh_s^x)ds+\int_{t_{n-1}}^{t_n} Dg_\eps (s,\Xh_s^x)\cdot dW_s\right)P_n\right]\right],
\e*
where $\E_{t_{j}}[\cdot]=\E[\cdot|\Fc_{t_{j}}]$. Notice that we can write $P_n=1+h^\frac12\beta_{n-1}\cdot\int_{t_{n-1}}^{t_n}dW_s +h^{-1}\alpha_{n-1}\cdot\int_{t_{n-1}}^{t_n}(W_s-W_{t_{n-1}})dW_s^\text{T}$. Thus, one can calculate $
\E_{t_{n-1}}\left[\left(\int_{t_{n-1}}^{t_n} \Lc^\Xh g_\eps (s,\Xh_s^x)ds+\int_{t_{n-1}}^{t_n} Dg_\eps (s,\Xh_s^x)\cdot dW_s\right)P_n\right]$ using  It\^o isometry and the fact that the expected value of stochastic integrals is zero:
\b*
&\E_{t_{n-1}}\left[\left(\int_{t_{n-1}}^{t_n} \Lc^\Xh g_\eps (s,\Xh_s^x)ds+\int_{t_{n-1}}^{t_n} Dg_\eps (s,\Xh_s^x)\cdot dW_s\right)P_n\right]\\
&=\E_{t_{n-1}}\left[\int_{t_{n-1}}^{t_n} \Lc^\Xh g_\eps (s,\Xh_s^x)ds+h^{\frac12}\beta_{n-1}\cdot \int_{t_{n-1}}^{t_n} Dg_\eps (s,\Xh_s^x)ds+h^{-1} \int_{t_{n-1}}^{t_n} \alpha_{n-1}Dg_\eps (s,\Xh_s^x)\cdot W_sds\right]
\e*
Because of \reff{g-eps-bnd}, the first two term in the above are bounded by $C(h+\frac{h}{\eps})$. The third term can be calculated by using \reff{hermit}
\b*
\E_{t_{n-1}}\left[\int_{t_{n-1}}^{t_n}  \alpha_{n-1}Dg_\eps (s,\Xh_s^x)\cdot W_sds\right]
&=&
\int_{t_{n-1}}^{t_n} \alpha_{n-1}\E_{t_{n-1}} [Dg_\eps (s,\Xh_s^x)\cdot W_s]ds
\e* 
By \reff{int-by-part}, we have $\E_{t_{n-1}} [Dg_\eps (s,\Xh_s^x)W_s]=s\E_{t_{n-1}} [D^2g_\eps (s,\Xh_s^x)]$ which is bounded by $C\frac{s}{\eps}$. Thus,
\b*
\left|\E_{t_{n-1}}\left[\int_{t_{n-1}}^{t_n}  \alpha_{n-1}Dg_\eps (s,\Xh_s^x)\cdot W_sds\right]\right|
&=&
\frac{C}{\eps}\int_{t_{n-1}}^{t_n}sds=\frac{Ch^2}{2\eps}.
\e* 
Therefore,
\b*
\E_{t_{n-1}}\left[\left(\int_{t_{n-1}}^{t_n} \Lc^\Xh g_\eps (s,\Xh_s^x)ds+\int_{t_{n-1}}^{t_n} Dg_\eps (s,\Xh_s^x)\cdot dW_s\right)P_n\right]\le 
Ch(1+ {\eps}^{-1}).
\e*
\no{\bf Step 5.}
Because $\{Q_j\}_{j=i}^n$ is a martingale and $\1_{\{t_{n} =\tauh\}}$ is $\Fc_{t_{n-1}}$-measurable, one can write the first term in the right hand side of \reff{bad-ass-term} as
\begin{equation*}
\begin{split}
&\E_{t_i,x}\left[\left(\int_{t_i}^{t_{n-1}} \Lc^\Xh g_\eps (s,\Xh_s^x)ds+\int_{t_i}^{t_{n-1}} Dg_\eps (s,\Xh_s^x)\cdot dW_s\right)Q_{n}\1_{\{t_{n} =\tauh\}}\right]\\
&=\E_{t_i,x}\left[\left(\int_{t_i}^{t_{n-1}} \Lc^\Xh g_\eps (s,\Xh_s^x)ds+\int_{t_i}^{t_{n-1}} Dg_\eps (s,\Xh_s^x)\cdot dW_s\right)\1_{\{t_{n} =\tauh\}}\E_{t_{n-1}}[Q_{n}]\right]\\
&=\E_{t_i,x}\left[\left(\int_{t_i}^{t_{n-1}} \Lc^\Xh g_\eps (s,\Xh_s^x)ds+\int_{t_i}^{t_{n-1}} Dg_\eps (s,\Xh_s^x)\cdot dW_s\right)Q_{n-1}\1_{\{t_{n} =\tauh\}}\right],
\end{split}
\end{equation*}
Thus, from \reff{summation-until-tau} we have
\begin{equation*}
\begin{split}
\Biggl|\E_{t_i,x}&\left[(g _\eps(\tauh,\Xh_\tauh^x)-g_\eps(t_i,x))\sum_{j=i}^{n-1}Q_{j+1}\1_{\{t_{j+1} =\tauh\}}\right]\Biggr|\le Ch(1+\frac{1}{\eps})\\
&+
\Biggl|\sum_{j=i}^{n-2}\E_{t_i,x}\left[\left(\int_{t_i}^\tauh \Lc^\Xh g_\eps (s,\Xh_s^x)ds+\int_{t_i}^\tauh Dg_\eps (s,\Xh_s^x)\cdot dW_s\right)Q_{j+1}\1_{\{t_{j+1} =\tauh\}}\right]\\
&+\E_{t_i,x}\left[\left(\int_{t_i}^{t_{n-1}} \Lc^\Xh g_\eps (s,\Xh_s^x)ds+\int_{t_i}^{t_{n-1}} Dg_\eps (s,\Xh_s^x)\cdot dW_s\right)Q_{n-1}\1_{\{t_{n} = \tauh\}}\right]\Biggr|\\
&=Ch(1+\frac{1}{\eps})+
\Biggl|\sum_{j=i}^{n-3}\E_{t_i,x}\left[\left(\int_{t_i}^\tauh \Lc^\Xh g_\eps (s,\Xh_s^x)ds+\int_{t_i}^\tauh Dg_\eps (s,\Xh_s^x)\cdot dW_s\right)Q_{j+1}\1_{\{t_{j+1} =\tauh\}}\right]\\
&+\E_{t_i,x}\left[\left(\int_{t_i}^{t_{n-1}} \Lc^\Xh g_\eps (s,\Xh_s^x)ds+\int_{t_i}^{t_{n-1}} Dg_\eps (s,\Xh_s^x)\cdot dW_s\right)Q_{n-1}\1_{\{t_{n-1} \le \tauh\}}\right]\Biggr|.
\end{split}
\end{equation*}
By  repeating the argument in {Step 3} and {Step 4} inductively over $k=n-1,\cdots,i+1$, one can write the above as
\begin{equation*}
\begin{split}
\Biggl|\E_{t_i,x}&\left[(g _\eps(\tauh,\Xh_\tauh^x)-g_\eps(t_i,x))\sum_{j=i}^{n-1}Q_{j+1}\1_{\{t_{j+1} =\tauh\}}\right]\Biggr|\le C(n-k)h(1+\frac{1}{\eps})\\
&+
\Biggl|\sum_{j=i}^{k-2}\E_{t_i,x}\left[\left(\int_{t_i}^\tauh \Lc^\Xh g_\eps (s,\Xh_s^x)ds+\int_{t_i}^\tauh Dg_\eps (s,\Xh_s^x)\cdot dW_s\right)Q_{j+1}\1_{\{t_{j+1} =\tauh\}}\right]\\
&+\E_{t_i,x}\left[\left(\int_{t_i}^{t_{k}} \Lc^\Xh g_\eps (s,\Xh_s^x)ds+\int_{t_i}^{t_{k}} Dg_\eps (s,\Xh_s^x)\cdot dW_s\right)Q_{k}\1_{\{t_{k} \le \tauh\}}\right]\Biggr|.
\end{split}
\end{equation*}
Specially for $k=i+1$ (the term containing $\sum_{j=i}^{k-2}$ disappears), we have 
\begin{equation*}
\begin{split}
\Biggl|\E_{t_i,x}&\left[(g _\eps(\tauh,\Xh_\tauh^x)-g_\eps(t_i,x))\sum_{j=i}^{n-1}Q_{j+1}\1_{\{t_{j+1} =\tauh\}}\right]\Biggr|\le C(n-i-1)h(1+\frac{1}{\eps})\\
&+
\E_{t_i,x}\left[\left(\int_{t_i}^{t_{i+1}} \Lc^\Xh g_\eps (s,\Xh_s^x)ds+\int_{t_i}^{t_{i+1}} Dg_\eps (s,\Xh_s^x)\cdot dW_s\right)Q_{i+1}\1_{\{t_{i+1} \le \tauh\}}\right]\Biggr|\\
&\le
C(n-i)h(1+\frac{1}{\eps})=C(T-t_i)(1+{\eps}^{-1}).
\end{split}
\end{equation*}
\no{\bf Step 6.} By using \reff{g-eps-bnd} and the bound found in Step 5 in \reff{vh-g}, one has 
\b*
|v^h(t_i,x)-g(t_i,x)|
&\le&
C(\eps+\frac{T-t_i}{\eps}+T-t_i).
\e*
By choosing $\eps=\sqrt{T-t_i}$, we conclude  that 
\b*
|v^h(t_i,x)-g(t_i,x)|\le C\sqrt{T-t_i}.
\e*
Then, the result follows from $x$-Lipschitz continuity and $t$-$\frac12$-H\"older continuity of $g$.} \ep

\no{\color{black}{\bf Acknowledgment.}
The authors are grateful to Xavier Warin and anonymous referees for their helpful comments and suggestions.}

\bibliographystyle{plain}	
\bibliography{bib0916}		

\begin{thebibliography}{10}

\bibitem{barlesjakobsen2}
G.~Barles and E.~R. Jakobsen.
\newblock On the convergence rate of approximation schemes for
  {H}amilton-{J}acobi-{B}ellman equations.
\newblock {\em M2AN Math. Model. Numer. Anal.}, 36(1):33--54, 2002.

\bibitem{barlesjakobsen1}
G.~Barles and E.~R. Jakobsen.
\newblock Error bounds for monotone approximation schemes for
  {H}amilton-{J}acobi-{B}ellman equations.
\newblock {\em SIAM J. Numer. Anal.}, 43(2):540--558 (electronic), 2005.

\bibitem{barlesjakobsen}
G.~Barles and E.~R. Jakobsen.
\newblock Error bounds for monotone approximation schemes for parabolic
  {H}amilton-{J}acobi-{B}ellman equations.
\newblock {\em Math. Comp.}, 76(260):1861--1893 (electronic), 2007.

\bibitem{barlessouganidis}
G.~Barles and P.~E. Souganidis.
\newblock Convergence of approximation schemes for fully nonlinear second order
  equations.
\newblock {\em Asymptotic Anal.}, 4(3):271--283, 1991.

\bibitem{2010arXiv1009.0932B}
Erhan Bayraktar and Yu-Jui Huang.
\newblock On the {M}ultidimensional {C}ontroller-and-{S}topper {G}ames.
\newblock {\em SIAM J. Control Optim.}, 51(2):1263--1297, 2013.

\bibitem{bmz}
J.~F. Bonnans, S.~Maroso, and H.~Zidani.
\newblock Error estimates for stochastic differential games: the adverse
  stopping case.
\newblock {\em IMA J. Numer. Anal.}, 26(1):188--212, 2006.

\bibitem{bc}
B.~Bouchard and J.-F. Chassagneux.
\newblock Discrete-time approximation for continuously and discretely reflected
  {BSDE}s.
\newblock {\em Stochastic Process. Appl.}, 118(12):2269--2293, 2008.

\bibitem{bt}
B.~Bouchard and N.~Touzi.
\newblock Discrete-time approximation and {M}onte-{C}arlo simulation of
  backward stochastic differential equations.
\newblock {\em Stochastic Process. Appl.}, 111(2):175--206, 2004.

\bibitem{bw}
B.~Bouchard and X.~Warin.
\newblock {Monte-Carlo valorisation of American options: facts and new
  algorithms to improve existing methods}.
\newblock {\em To appear in Numerical Methods in Finance , Springer Proceedings
  in Mathematics, ed. R. Carmona, P. Del Moral, P. Hu and N. Oudjane}, 2011.

\bibitem{caffarelli-souganidis}
Luis~A Caffarelli and Panagiotis~E Souganidis.
\newblock A rate of convergence for monotone finite difference approximations
  to fully nonlinear, uniformly elliptic pdes.
\newblock {\em Communications on Pure and Applied Mathematics}, 61(1):1--17,
  2008.

\bibitem{cil}
M.~G. Crandall, H.~Ishii, and P.-L. Lions.
\newblock User's guide to viscosity solutions of second order partial
  differential equations.
\newblock {\em Bull. Amer. Math. Soc. (N.S.)}, 27(1):1--67, 1992.

\bibitem{fahim}
A.~Fahim.
\newblock Convergence of a {M}onte {C}arlo method for fully non-linear elliptic
  and parabolic {PDE}s in some general domains.
\newblock Preprint, Sep 2011.

\bibitem{ftw}
A.~Fahim, N.~Touzi, and X.~Warin.
\newblock A probabilistic numerical method for fully non--linear parabolic
  pdes.
\newblock {\em Annals of Applied Probability}, 21(4):1322--1364, 2011.

\bibitem{gzz}
Wenjie Guo, Jianfeng Zhang, and Jia Zhuo.
\newblock A monotone scheme for high dimensional fully nonlinear pdes.
\newblock {\em arXiv preprint arXiv:1212.0466}, 2012.

\bibitem{jakobsen}
E.~R. Jakobsen.
\newblock On the rate of convergence of approximation schemes for {B}ellman
  equations associated with optimal stopping time problems.
\newblock {\em Math. Models Methods Appl. Sci.}, 13(5):613--644, 2003.

\bibitem{MR1809521}
I.~Karatzas and S.~G. Kou.
\newblock Hedging {A}merican contingent claims with constrained portfolios.
\newblock {\em Finance Stoch.}, 2(3):215--258, 1998.

\bibitem{krylov}
N.~V. Krylov.
\newblock On the rate of convergence of finite-difference approximations for
  {B}ellman's equations.
\newblock {\em Algebra i Analiz}, 9(3):245--256, 1997.

\bibitem{krylov0}
N.~V. Krylov.
\newblock Approximating value functions for controlled degenerate diffusion
  processes by using piece-wise constant policies.
\newblock {\em Electron. J. Probab.}, 4:no. 2, 19 pp. (electronic), 1999.

\bibitem{krylov2}
N.~V. Krylov.
\newblock On the rate of convergence of finite-difference approximations for
  {B}ellman's equations with variable coefficients.
\newblock {\em Probab. Theory Related Fields}, 117(1):1--16, 2000.

\bibitem{krylov1}
N.~V. Krylov.
\newblock The rate of convergence of finite-difference approximations for
  {B}ellman equations with {L}ipschitz coefficients.
\newblock {\em Appl. Math. Optim.}, 52(3):365--399, 2005.

\bibitem{longsch}
F.A. Longstaff and E.~S. Schwartz.
\newblock Valuing {A}merican options by simulation: a simple least-squares
  approach.
\newblock {\em Review of Financial Studies}, 14(1):113--147, 2001.

\bibitem{mazhang}
J.~Ma and J.~Zhang.
\newblock Representations and regularities for solutions to {BSDE}s with
  reflections.
\newblock {\em Stochastic Process. Appl.}, 115(4):539--569, 2005.

\bibitem{oberman}
A.~M. Oberman.
\newblock Convergent difference schemes for degenerate elliptic and parabolic
  equations: {H}amilton-{J}acobi equations and free boundary problems.
\newblock {\em SIAM J. Numer. Anal.}, 44(2):879--895 (electronic), 2006.

\bibitem{obermanthaleia}
A.~M. Oberman and T.~Zariphopoulou.
\newblock Pricing early exercise contracts in incomplete markets.
\newblock {\em Computational Management Science}, 1(1):75--107, 2003.

\bibitem{zhang}
J.~Zhang.
\newblock A numerical scheme for {BSDE}s.
\newblock {\em Ann. Appl. Probab.}, 14(1):459--488, 2004.

\end{thebibliography}

\end{document}